\RequirePackage{fix-cm}
\documentclass{siamart1116}
\usepackage{amsmath}
\usepackage{amssymb}
\usepackage{mathtools}

\usepackage{pgfplots, pgfplotstable, booktabs, colortbl, array}
\usepackage{tikz}
\pgfplotsset{compat=newest}
\usepackage{pgffor}
\usepackage{xcolor}

\DeclareMathOperator*{\argmin}{arg\,min}

\renewcommand{\Re}{\operatorname{Re}}
\renewcommand{\Im}{\operatorname{Im}}

\DeclarePairedDelimiter\floor{\lfloor}{\rfloor}

\usepackage{tikz}
\usetikzlibrary{positioning,arrows}

\newcommand{\maketable}[1]{%
\pgfplotstabletypeset[
every head row/.style={before row=\toprule,after row=\midrule},
every last row/.style={after row=\bottomrule},
columns={leftcol,0,1,2},
create on use/leftcol/.style={create col/set list={DB,DBp,IN,{P-$(1,0)$},{P-$(4,4)$},{P-$(8,8)$},{Z-$(1,0)$},{Z-$(4,4)$},{Z-$(8,8)$},sqrtm}},
columns/leftcol/.style={string type,column type/.add={|}{|},column name={Method}},
columns/0/.style={fixed,column type/.add={}{|},column name={$k$}},
columns/1/.style={sci,sci e,sci zerofill,precision=1,verbatim,column type/.add={}{|},column name={Err.}},
columns/2/.style={sci,sci e,sci zerofill,precision=1,verbatim,column type/.add={}{|},column name={Res.}}
]
{#1}
}

\begin{document}

\title{Zolotarev Iterations for the Matrix Square Root\thanks{Submitted to the editors March 31, 2018.}}

\author{Evan S. Gawlik\thanks{Department of Mathematics, University of California, San Diego
    (\email{egawlik@ucsd.edu})}}

\date{}

\headers{Zolotarev iterations for the matrix square root}{E. S. Gawlik}

\maketitle

\begin{abstract}
We construct a family of iterations for computing the principal square root of a square matrix $A$ using Zolotarev's rational minimax approximants of the square root function.  We show that these rational functions obey a recursion, allowing one to iteratively generate optimal rational approximants of $\sqrt{z}$ of high degree using compositions and products of low-degree rational functions.  The corresponding iterations for the matrix square root converge to $A^{1/2}$ for any input matrix $A$ having no nonpositive real eigenvalues.  In special limiting cases, these iterations reduce to known iterations for the matrix square root: the lowest-order version is an optimally scaled Newton iteration, and for certain parameter choices, the principal family of Pad\'e iterations is recovered.  Theoretical results and numerical experiments indicate that the iterations perform especially well on matrices having eigenvalues with widely varying magnitudes.
\end{abstract}

\section{Introduction} \label{sec:intro}

A well-known method for computing the square root of an $n \times n$ matrix $A$ with no nonpositive real eigenvalues is the Newton iteration~\cite{higham1986newton}
\begin{equation} \label{newton}
X_{k+1} = \frac{1}{2}(X_k + X_k^{-1}A), \quad X_0 = A.
\end{equation}
In exact arithmetic, the matrix $X_k$ converges quadratically to $A^{1/2}$, the  principal square root of $A$~\cite[Theorem 6.9]{higham2008functions}.  (In floating point arithmetic, mathematically equivalent reformulations of~(\ref{newton}), such as the Denman-Beavers iteration~\cite{denman1976matrix}, are preferred for stability reasons~\cite[Section 6.4]{higham2008functions}.)

If $A$ is diagonalizable, then each eigenvalue $\lambda^{(i)}_k$ of $X_k$, $i=1,2,\dots,n$, obeys a recursion of the form
\[
\lambda_{k+1}^{(i)} = \frac{1}{2}\left(\lambda_k^{(i)} + \frac{\lambda_0^{(i)}}{\lambda_k^{(i)}}\right),
\]
which is the Newton iteration for computing a root of $z^2-\lambda_0^{(i)}=0$.
One can thus think of~(\ref{newton}) as an iteration that, in the limit as $k \rightarrow \infty$, implicitly maps a collection of scalars $\lambda_0^{(i)}$, $i=1,2,\dots,n$, to $\sqrt{\lambda_0^{(i)}}$ for each $i$.  In order for each scalar to converge rapidly, it is necessary that the rational function $f_k(z)$ defined recursively by
\begin{equation} \label{rationalnewton}
f_{k+1}(z) = \frac{1}{2}\left(f_k(z) + \frac{z}{f_k(z)}\right), \quad f_0(z) = z
\end{equation}
converges rapidly to $f(z) = \sqrt{z}$ on the set $\bigcup_{i=1}^n \{\lambda_0^{(i)}\} \subset \mathbb{C}$.

To generalize and improve the Newton iteration, it is natural to study other recursive constructions of rational functions, with the aim of approximating $\sqrt{z}$ on a subset $S \subset \mathbb{C}$ containing the spectrum of $A$.  Of particular interest are rational functions that minimize the maximum relative error
\begin{equation} \label{maxerr}
\max_{z \in S} \left|(r(z)-\sqrt{z})/\sqrt{z}\right|
\end{equation}
among all rational functions $r(z)$ of a given type $(m,\ell)$.  By type $(m,\ell)$, we mean $r(z)=p(z)/q(z)$ is a ratio of polynomials $p$ and $q$ of degree at most $m$ and $\ell$, respectively. 
We denote the set of rational functions of type $(m,\ell)$ by $\mathcal{R}_{m,\ell}$.

On a positive real interval $S$, explicit formulas for the minimizers $r \in \mathcal{R}_{m,m-1}$ and $r \in \mathcal{R}_{m,m}$ of~(\ref{maxerr}) are known for each $m$.  The formulas, derived by Zolotarev~\cite{zolotarev1877applications}, are summarized in Section~\ref{sec:background}.  We show in this paper that, remarkably, the minimizers obey a recursion analogous to~(\ref{rationalnewton}).  This fact is intimately connected to (and indeed follows from) an analogous recursion for rational minimax approximations of the function $\mathrm{sign}(z) = z/\sqrt{z^2}$ recently discovered by Nakatsukasa and Freund~\cite{nakatsukasa2016computing}.  

The lowest order version of the recursion for square root approximants has been known for several decades~\cite{rutishauser1963betrachtungen,ninomiya1970best}~\cite[Section V.5.C]{braess1986nonlinear}. Beckermann~\cite{beckermann2013optimally} recently studied its application to matrices.
In this paper, we generalize these ideas by constructing a family of iterations for computing the matrix square root, one for each pair of integers $(m,\ell)$ with $\ell \in \{m-1,m\}$.  
We prove that these \emph{Zolotarev iterations} are stable and globally convergent with order of convergence $m+\ell+1$.  By writing Zolotarev's rational functions in partial fraction form, the resulting algorithms are highly parallelizable.  Numerical examples demonstrate that the iterations exhibit good forward stability. 

The Zolotarev iterations for the matrix square bear several similarities to the Pad\'{e} iterations studied in~\cite[pp. 231-233]{higham1997stable},~\cite[Section 6]{higham2005functions}, and~\cite[Section 6.7]{higham2008functions}.  In fact, the Pad\'e iterations can be viewed as a limiting case of the Zolotarev iterations; see Proposition~\ref{prop:pade}.  One of the messages we hope to convey in this paper is that the Zolotarev iterations are often preferable to the Pad\'e iterations when the eigenvalues of $A$ have widely varying magnitudes.  Roughly, this can be understood by noting that the Pad\'e approximants of $\sqrt{z}$ are designed to be good approximations of $\sqrt{z}$ near a point, whereas Zolotarev's minimax approximants are designed to be good approximations of $\sqrt{z}$ over an entire interval.  For more details, particularly with regards to how these arguments carry over to the complex plane, see Section~\ref{sec:scalar}.

This paper builds upon a stream of research that, in recent years, has sparked renewed interest in the applications of Zolotarev's work on rational approximation to numerical linear algebra.  These applications include algorithms for the SVD, the symmetric eigendecomposition, and the polar decomposition~\cite{nakatsukasa2016computing}; algorithms for the CS decomposition~\cite{gawlik2018backward}; bounds on the singular values of matrices with displacement structure~\cite{beckermann2017singular}; computation of spectral projectors~\cite{kressner2017fast,guttel2015zolotarev,li2017spectrum}; and the selection of optimal parameters for the alternating direction implicit (ADI) method~\cite{wachspress2013adi,bailly2000optimal}.  Zolotarev's functions have even been used to compute the matrix square root~\cite{hale2008computing}, however, there is an important distinction between that work and ours: In~\cite{hale2008computing}, Zolotarev's functions are not used as the basis of an iterative method.  Rather, a rational function of $A$ is evaluated once and for all to approximate $A^{1/2}$.  As we argue below, recursive constructions of Zolotarev's functions offer significant advantages over this strategy.  

This paper is organized as follows.  In Section~\ref{sec:results}, we state our main results without proof.  In Section~\ref{sec:proofs}, we prove these results.  In Section~\ref{sec:practical}, we discuss the implementation of the Zolotarev iterations and how they compare with other known iterations.  In Section~\ref{sec:numerical}, we evaluate the performance of the Zolotarev iterations on numerical examples.

\section{Statement of Results} \label{sec:results}

In this section, we state our main results and discuss some of their implications.  Proofs are presented in Section~\ref{sec:proofs}.

\paragraph{Recursion for rational approximations of $\sqrt{z}$}
We begin by introducing a recursion satisfied by Zolotarev's best rational approximants of the square root function.  For each $m \in \mathbb{N}$, $\ell \in \{m-1,m\}$, and $\alpha \in (0,1)$, let $r_{m,\ell}(z,\alpha)$ denote the rational function of type $(m,\ell)$ that minimizes~(\ref{maxerr}) on $S=[\alpha^2,1]$.  Let $\hat{r}_{m,\ell}(z,\alpha)$ be the unique scalar multiple of $r_{m,\ell}(z,\alpha)$ with the property that 
\[
\min_{z \in [\alpha^2,1]} (\hat{r}_{m,\ell}(z,\alpha)-\sqrt{z})/\sqrt{z}  = 0.
\]
The following theorem, which is closely related to~\cite[Corollary 4]{nakatsukasa2016computing} and includes~\cite[Lemma 1]{beckermann2013optimally} as a special case, will be proved in Section~\ref{sec:proofs}.
\begin{theorem} \label{thm:rationalzolo}
Let $m \in \mathbb{N}$ and $\alpha \in (0,1)$.  Define $f_k(z)$ recursively by
\begin{align} 
f_{k+1}(z) &= f_k(z) \hat{r}_{m,m-1}\left( \frac{z}{f_k(z)^2}, \alpha_k \right), & f_0(z) &= 1, \label{rationalzoloA1} \\
\alpha_{k+1} &= \frac{\alpha_k}{\hat{r}_{m,m-1}(\alpha_k^2,\alpha_k)}, & \alpha_0 &= \alpha. \label{rationalzoloA2}
\end{align}
Then, for every $k \ge 1$,
\[
f_k(z) = \hat{r}_{q,q-1}(z,\alpha) = \frac{1+\alpha_k}{2\alpha_k} r_{q,q-1}(z,\alpha), \quad q = \frac{1}{2}(2m)^k.
\]
If instead
\begin{align}
f_{k+1}(z) &= f_k(z) \hat{r}_{m,m}\left( \frac{z}{f_k(z)^2}, \alpha_k \right), & f_0(z) &= 1, \label{rationalzoloB1} \\
\alpha_{k+1} &= \frac{\alpha_k}{\hat{r}_{m,m}(\alpha_k^2,\alpha_k)}, & \alpha_0 &= \alpha, \label{rationalzoloB2}
\end{align}
then, for every $k \ge 1$,
\[
f_k(z) = \hat{r}_{q,q}(z,\alpha) = \frac{1+\alpha_k}{2\alpha_k} r_{q,q}(z,\alpha), \quad q = \frac{1}{2}((2m+1)^k-1).
\]
\end{theorem}

The remarkable nature of these recursions is worth emphasizing with an example.  When $m=7$, three iterations of~(\ref{rationalzoloA1}-\ref{rationalzoloA2}) generate (up to rescaling) the best rational approximation of $\sqrt{z}$ of type $(1372,1371)$ on the interval $[\alpha^2,1]$.
Not only is this an efficient way of computing $r_{1372,1371}(z,\alpha)$, but it also defies intuition that an iteration involving so few parameters could deliver the solution to an optimization problem (the minimization of~(\ref{maxerr}) over $\mathcal{R}_{1372,1371}$) with thousands of degrees of freedom. 

\paragraph{Zolotarev iterations for the matrix square root} Theorem~\ref{thm:rationalzolo} leads to a family of iterations for computing the square root of an $n \times n$ matrix $A$, namely,
\begin{align}
X_{k+1} &= X_k \hat{r}_{m,\ell}(X_k^{-2}A, \alpha_k), & X_0 &= I, \label{zolo1} \\
\alpha_{k+1} &= \frac{\alpha_k}{\hat{r}_{m,\ell}(\alpha_k^2,\alpha_k)}, & \alpha_0 &= \alpha, \label{zolo2}
\end{align}
where $m$ is a positive integer and $\ell \in \{m-1,m\}$.
We will refer to each of these iterations as a \emph{Zolotarev iteration} of type $(m,\ell)$.
(Like the Newton iteration, these iterations are ill-suited for numerical implementation in their present form, but a reformulation renders them numerically stable; see the end of this section).
A priori, these iterations would appear to be suitable only for Hermitian positive definite matrices (or, more generally, diagonalizable matrices with positive real eigenvalues) that have been scaled so that their eigenvalues lie in the interval $[\alpha^2,1]$, 
but in fact they converge for any $A \in \mathbb{C}^{n \times n}$ with no nonpositive real eigenvalues.  This is made precise in the forthcoming theorem, which is a generalization of~\cite[Theorem 4]{beckermann2013optimally} and is related to~\cite[Theorem 4.1]{hale2008computing}.

To state the theorem, we introduce some notation, following~\cite{beckermann2013optimally}. A compact set $S \subseteq \mathbb{C}$ is called $L$-spectral for $A \in \mathbb{C}^{n \times n}$ if
\[
\|f(A)\|_2 \le L \sup_{z \in S} |f(z)|
\]
for every function $f$ analytic in $S$~\cite[Chapter 37]{hogben2016handbook}.  For instance, the spectrum of $A$ is $1$-spectral for every normal matrix $A$, and the closure of the pseudospectrum $\Lambda_{\epsilon}(A) = \{z \in \mathbb{C} \mid \|(A-zI)^{-1}\|_2 > 1/\epsilon\}$ is $C_\epsilon$-spectral with $C_\epsilon = \mathrm{length}(\partial\Lambda_\epsilon(A))/(2\pi\epsilon)$ for every $A$~\cite[Fact 23.3.5]{hogben2016handbook}.  

For each $\alpha \in (0,1)$, define
\begin{equation} \label{phi}
\varphi(z,\alpha) = \exp\left( \frac{ \pi \mathrm{sn}^{-1}(\sqrt{z}/\alpha; \alpha) }{ K(\alpha') } \right),
\end{equation}
where $\mathrm{sn}(\cdot;\alpha)$, $\mathrm{cn}(\cdot;\alpha)$, and $\mathrm{dn}(\cdot;\alpha)$ denote Jacobi's elliptic functions with modulus $\alpha$, $K(\alpha) = \int_0^{\pi/2} (1-\alpha^2 \sin^2\theta)^{-1/2} \, d\theta$ is the complete elliptic integral of the first kind, and $\alpha' = \sqrt{1-\alpha^2}$ is the complementary modulus to $\alpha$.  Note that the function $\varphi(z,\alpha)$ supplies a conformal map from $\mathbb{C} \setminus ((-\infty,0] \cup [\alpha^2,1])$ to the annulus $\{z \in \mathbb{C} : 1 < |z| < \rho(\alpha)\}$~\cite[pp. 138-140]{akhiezer1990elements}, where
\begin{equation} \label{rho}
\rho(\alpha) = \exp\left( \frac{\pi K(\alpha)}{K(\alpha')} \right).
\end{equation}

\begin{theorem} \label{thm:convergence}
Let $A \in \mathbb{C}^{n \times n}$ have no nonpositive real eigenvalues.  
Suppose that $S \subseteq \mathbb{C} \setminus (-\infty,0]$ is $L$-spectral for $A$.  
Let $m \in \mathbb{N}$, $\ell \in \{m-1,m\}$, $\alpha \in (0,1)$, and $\gamma = \inf_{z \in S} |\varphi(z,\alpha)|$.  For every $k \ge 1$ such that $\max\{2\gamma^{-2(m+\ell+1)^k},4\rho(\alpha)^{-2(m+\ell+1)^k}\} < 1$, the matrix $X_k$ defined by~(\ref{zolo1}-\ref{zolo2}) satisfies
\begin{align}
\left\| \left( \frac{2\alpha_k}{1+\alpha_k} \right) X_k A^{-1/2} - I \right\|_2 &\le 4 L \gamma^{-(m+\ell+1)^k} + O\left( \gamma^{-2(m+\ell+1)^k} \right). \label{matest}
\end{align}
If $S \subseteq [\alpha^2,1]$, then~(\ref{matest}) holds with $O\left( \gamma^{-2(m+\ell+1)^k} \right)$ replaced by zero, $\gamma = \rho(\alpha)$, and $k \ge 1$.
\end{theorem}

\begin{corollary} \label{cor:hermitian}
Let $A \in \mathbb{C}^{n \times n}$ be Hermitian positive definite.  If the eigenvalues of $A$ lie in the interval $[\alpha^2,1]$, then
\[
\left\| \left( \frac{2\alpha_k}{1+\alpha_k} \right) X_k A^{-1/2} - I \right\|_2 \le 4 \rho(\alpha)^{-(m+\ell+1)^k} 
\]
for every $k \ge 1$.
\end{corollary}

Note that the error estimates above imply estimates for the relative error in the computed square root $\widetilde{X}_k := 2\alpha_k X_k / (1+\alpha_k)$, since
\[
\frac{\|\widetilde{X}_k-A^{1/2}\|_2}{\|A^{1/2}\|_2} = \frac{\|(\widetilde{X}_k A^{-1/2} - I) A^{1/2}\|_2}{\|A^{1/2}\|_2} \le  \|\widetilde{X}_k A^{-1/2} - I\|_2.
\]

\paragraph{Connections with existing iterations}
It is instructive to examine the lowest order realization of the iteration~(\ref{zolo1}-\ref{zolo2}).  When $(m,\ell)=(1,0)$, one checks (using either elementary calculations or the explicit formulas in Section~\ref{sec:background}) that
\[
\hat{r}_{1,0}(z,\alpha) = \frac{1}{2}(\alpha^{1/2} + \alpha^{-1/2}z),
\]
so that the iteration~(\ref{zolo1}-\ref{zolo2}) reduces to
\begin{align*}
X_{k+1} &= \frac{1}{2}(\alpha_k^{1/2} X_k + \alpha_k^{-1/2} X_k^{-1} A), & X_0 &= I, \\
\alpha_{k+1} &= \frac{2}{\alpha_k^{1/2}+\alpha_k^{-1/2}}, & \alpha_0 &= \alpha.
\end{align*}
Equivalently, in terms of $\mu_k := \alpha_k^{1/2}$,
\begin{align}
X_{k+1} &= \frac{1}{2}(\mu_k X_k + \mu_k^{-1} X_k^{-1} A), & X_0 &= I, \label{scaledNewton1} \\
\mu_{k+1} &= \sqrt{\frac{2}{\mu_k+\mu_k^{-1}}}, & \mu_0 &= \alpha^{1/2}. \label{scaledNewton2}
\end{align}
This is precisely the scaled Newton iteration with a scaling heuristic studied in~\cite{beckermann2013optimally}.  (In~\cite{beckermann2013optimally}, starting values $X_0=A$ and $\mu_0=\alpha^{-1/2}$ are used, but it easy to check that this generates the same sequences $\{X_k\}_{k=1}^\infty$ and $\{\mu_k\}_{k=1}^\infty$ as~(\ref{scaledNewton1}-\ref{scaledNewton2}).)   This iteration has its roots in early work on rational approximation of the square root~\cite{rutishauser1963betrachtungen,ninomiya1970best}, and it is closely linked to the scaled Newton iteration for the polar decomposition introduced in~\cite{byers2008new}.  As with the unscaled Newton iteration, reformulating~(\ref{scaledNewton1}-\ref{scaledNewton2}) (e.g., as a scaled Denman-Beavers iteration) is necessary to ensure its numerical stability.

Another class of known iterations for the matrix square root is recovered if one examines the limit as $\alpha \uparrow 1$.
Below, we say that a family of functions $\{r(\cdot,\alpha) \in \mathcal{R}_{m,\ell} : \alpha \in (0,1)\}$ converges \emph{coefficient-wise} to a function $p \in \mathcal{R}_{m,\ell}$ as $\alpha \uparrow 1$ if the coefficients of the polynomials in the numerator and denominator of $r(z,\alpha)$, appropriately normalized, approach the corresponding coefficients in $p(z)$ as $\alpha \uparrow 1$.
\begin{proposition} \label{prop:pade}
Let $m \in \mathbb{N}$ and $\ell \in \{m-1,m\}$.  As $\alpha \uparrow 1$, $\hat{r}_{m,\ell}(z,\alpha)$ converges coefficient-wise to $p_{m,\ell}(z)$, the type $(m,\ell)$ Pad\'{e} approximant of $\sqrt{z}$ at $z=1$. 
\end{proposition}

Since $p_{m,\ell}(1)=1$, the iteration~(\ref{zolo1}-\ref{zolo2}) formally reduces to
\begin{equation} \label{pade}
X_{k+1} = X_k p_{m,\ell}(X_k^{-2}A), \quad X_0 = A
\end{equation}  
as $\alpha \uparrow 1$.
To relate this to an existing iteration from the literature, define $Y_k = X_k^{-1}A$ and $Z_k = X_k^{-1}$.  Then, using the mutual commutativity of $X_k$, $Y_k$, $Z_k$, and $A$, we arrive at the iteration
\begin{align} 
Y_{k+1} &= Y_k q_{\ell,m}( Z_k Y_k), & Y_0 &= A, \label{padecoupled1} \\
Z_{k+1} &= q_{\ell,m}( Z_k Y_k ) Z_k, & Z_0 &= I, \label{padecoupled2}
\end{align} 
where $q_{\ell,m}(z) = p_{m,\ell}(z)^{-1}$.
Since $q_{\ell,m}(z)$ is the type $(\ell,m)$ Pad\'e approximant of $z^{-1/2}$ at $z=1$, this iteration is precisely the Pad\'e iteration studied in~\cite[Section 6]{higham2005functions},~\cite[p. 232]{higham1997stable}, and~\cite[Section 6.7]{higham2008functions}.  There, it is shown that $Y_k \rightarrow A^{1/2}$ and $Z_k \rightarrow A^{-1/2}$ with order of convergence $m+\ell+1$ for any $A$ with no nonpositive real eigenvalues.  Moreover, the iteration~(\ref{padecoupled1}-\ref{padecoupled2}) is stable~\cite[Theorem 6.12]{higham2008functions}.

\paragraph{Stable reformulation of the Zolotarev iterations}
In view of the well-established stability theory for iterations of the form~(\ref{padecoupled1}-\ref{padecoupled2}), we will focus in this paper on the following reformulation of the Zolotarev iteration~(\ref{zolo1}-\ref{zolo2}):
\begin{align} 
Y_{k+1} &= Y_k h_{\ell,m}( Z_k Y_k, \alpha_k), & Y_0 &= A, \label{zolocoupled1} \\
Z_{k+1} &= h_{\ell,m}( Z_k Y_k, \alpha_k ) Z_k, & Z_0 &= I, \label{zolocoupled2} \\
\alpha_{k+1} &= \alpha_k h_{\ell,m}(\alpha_k^2,\alpha_k), & \alpha_0 &= \alpha, \label{zolocoupled3}
\end{align}
where $h_{\ell,m}(z,\alpha) = \hat{r}_{m,\ell}(z,\alpha)^{-1}$.  In exact arithmetic, $Y_k$ and $Z_k$ are related to $X_k$ from~(\ref{zolo1}-\ref{zolo2}) via $Y_k = X_k^{-1}A$, $Z_k = X_k^{-1}$.  The following theorem summarizes the properties of this iteration.

\begin{theorem} \label{thm:stability}
Let $m \in \mathbb{N}$, $\ell \in \{m-1,m\}$, and $\alpha \in (0,1)$.  For any $A \in \mathbb{C}^{n \times n}$ with no nonpositive real eigenvalues, the iteration~(\ref{zolocoupled1}-\ref{zolocoupled3}) is stable, and $Y_k \rightarrow A^{1/2}$, $Z_k \rightarrow A^{-1/2}$, and $\alpha_k \rightarrow 1$ with order of convergence $m+\ell+1$.
\end{theorem}

Note that although Theorem~\ref{thm:stability} places no restrictions on the spectral radius of $A$ nor the choice of $\alpha \in (0,1)$, it should be clear that it is preferable to scale $A$ so that its spectral radius is 1 (or approximately 1), and set $\alpha = \sqrt{|\lambda_{\mathrm{min}}/\lambda_{\mathrm{max}}|}$ (or an estimate thereof), where $\lambda_{\mathrm{max}}$ and $\lambda_{\mathrm{min}}$ are the eigenvalues of $A$ with the largest and smallest magnitudes, respectively.  See Section~\ref{sec:scalar} for more details.

\section{Proofs} \label{sec:proofs}

In this section, we present proofs of Theorem~\ref{thm:rationalzolo}, Theorem~\ref{thm:convergence}, Proposition~\ref{prop:pade}, and Theorem~\ref{thm:stability}.

\subsection{Background} \label{sec:background}

We begin by reviewing a few facts from the theory of rational minimax approximation.  For a thorough presentation of this material, see, for example,~\cite[Chapter II]{akhiezer1956theory} and~\cite[Chapter 9]{akhiezer1990elements}.

\paragraph{Rational minimax approximation} 

Let $S = [a,b]$ be a finite interval.  A function $g(z)$ is said to \emph{equioscillate} between $N$ extreme points on $S$ if there exist $N$ points $z_1 < z_2 < \cdots < z_N$ in $S$ at which
\[
g(z_j) = \sigma (-1)^j \max_{z \in S} |g(z)|, \quad j=1,2,\dots,N,
\]
for some $\sigma \in \{-1,1\}$.  

Let $f$ and $w$ be continuous, real-valued functions on $S$ with $w>0$ on $S$.  Consider the problem of finding a rational function $r \in \mathcal{R}_{p,q}$ that minimizes
\[
\max_{z \in S} |(r(z)-f(z))w(z)|
\]
among all rational functions of type $(p,q)$.  It is well-known that this problem admits a unique solution $r^*$~\cite[p. 55]{akhiezer1956theory}.  Furthermore, the following are sufficient conditions guaranteeing optimality: If $r \in \mathcal{R}_{p,q}$ has the property that $(r(z)-f(z))w(z)$ equioscillates between $p+q+2$ extreme points on $S$, then $r=r^*$~\cite[p. 55]{akhiezer1956theory}.  (If $S$ is a union of two disjoint intervals, then this statement holds with $p+q+2$ replaced by $p+q+3$~\cite[Lemma 2]{nakatsukasa2016computing}.)

\paragraph{Rational approximation of the sign function}
Our analysis will make use of a connection between rational minimax approximants of $\sqrt{z}$ and rational minimax approximants of the function $\mathrm{sign}(z) = z/\sqrt{z^2}$.
For $\alpha \in (0,1)$, define
\[
s_{p,q}(z,\alpha) = \argmin_{s \in \mathcal{R}_{p,q}} \max_{z \in [-1,-\alpha] \cup [\alpha,1]} |s(z)-\mathrm{sign}(z)|
\]
and
\begin{equation} \label{errsign}
E_{p,q}(\alpha) = \max_{z \in [-1,-\alpha] \cup [\alpha,1]} |s_{p,q}(z,\alpha)-\mathrm{sign}(z)|.
\end{equation}
We will be primarily interested in the functions $s_{p,q}$ with $p \in \{2m-1,2m+1\}$ and $q=2m$, for which explicit formulas are known thanks to the seminal work of Zolotarev~\cite{zolotarev1877applications}. Namely, for $\ell \in \{m-1,m\}$, we have~\cite[p. 286]{akhiezer1956theory}
\[
s_{2\ell+1,2m}(z,\alpha) = M(\alpha) z \frac{\prod_{j=1}^{\ell} (z^2+c_{2j}(\alpha))}{\prod_{j=1}^m (z^2+c_{2j-1}(\alpha))},
\]
where
\begin{align}
c_j(\alpha) &= \alpha^2 \frac{\mathrm{sn}^2\left( \frac{jK(\alpha')}{m+\ell+1}; \alpha' \right)}{\mathrm{cn}^2 \left( \frac{jK(\alpha')}{m+\ell+1}; \alpha' \right)}, \label{cj}
\end{align}
and $M(\alpha)$ is a scalar uniquely defined by the condition that
\[
\min_{z \in [\alpha,1]} \left( s_{2\ell+1,2m}(z,\alpha) - 1 \right) = -\max_{z \in [\alpha,1]} \left( s_{2\ell+1,2m}(z,\alpha) - 1 \right).
\]

For $\ell \in \{m-1,m\}$, we denote
\[
\varepsilon_{m,\ell}(\alpha) = E_{2\ell+1,2m}(\alpha),
\]
and we use the abbreviation $\varepsilon := \varepsilon_{m,\ell}(\alpha)$ whenever there is no danger of confusion.  For each $\ell \in \{m-1,m\}$, it can be shown that
on the interval $[\alpha,1]$, $s_{2\ell+1,2m}(z,\alpha)$ takes values in $[1-\varepsilon,1+\varepsilon]$ and achieves its extremal values at exactly $m+\ell+2$ points $\alpha = z_0 < z_1 < \dots < z_{m+\ell+1} = 1$~\cite[p. 286]{akhiezer1956theory}:
\begin{equation} \label{equipts}
s_{2\ell+1,2m}(z_j,\alpha) - 1 = (-1)^{j+1} \varepsilon, \quad z_j = \alpha \bigg/ \mathrm{dn}\left( \frac{jK(\alpha')}{m+\ell+1}; \alpha' \right), \, j=0,1,\dots,m+\ell+1.
\end{equation} 
Since $s_{2\ell+1,2m}(z,\alpha)$ is odd, it follows that $s_{2\ell+1,2m}(z)-\mathrm{sign}(z)$ equioscillates between $2m+2\ell+4$ extreme points on $[-1,-\alpha] \cup [\alpha,1]$, confirming its optimality.

An important role in what follows will be played by the scaled function
\begin{equation} \label{shat}
\hat{s}_{2\ell+1,2m}(z,\alpha) = \frac{1}{1+\varepsilon_{m,\ell}(\alpha)} s_{2\ell+1,2m}(z,\alpha),
\end{equation}
which has the property that
\[
\max_{z \in [-1,-\alpha] \cup [\alpha,1]}
|\hat{s}_{2\ell+1,2m}(z,\alpha)| = 1.
\]

\paragraph{Rational approximation of the square root function}
For each $m$ and each $\ell \in \{m-1,m\}$, define
\[
r_{m,\ell}(z,\alpha) = (1-\varepsilon_{m,\ell}(\alpha)) (1+\varepsilon_{m,\ell}(\alpha)) \frac{\sqrt{z}}{s_{2\ell+1,2m}(\sqrt{z},\alpha)}
\]
and
\begin{equation} \label{rhat}
\hat{r}_{m,\ell}(z,\alpha) = \frac{1}{1-\varepsilon_{m,\ell}(\alpha)} r_{m,\ell}(z,\alpha).
\end{equation}
We claim that these definitions are consistent with those in Section~\ref{sec:intro}.  That is, $r_{m,\ell}(z,\alpha)$ minimizes
\[
\max_{z \in [\alpha^2,1]} |r(z)/\sqrt{z}-1|
\]
among all $r \in \mathcal{R}_{m,\ell}$, and $\hat{r}_{m,\ell}(z,\alpha)$ is scaled in such a way that
\[
\min_{z \in [\alpha^2,1]} (\hat{r}_{m,\ell}(z,\alpha)/\sqrt{z}-1) = 0.
\]
Indeed, it is easy to see from the properties of $s_{2\ell+1,2m}(z,\alpha)$ that (denoting $\varepsilon := \varepsilon_{m,\ell}(\alpha)$):
\begin{enumerate}
\item $r_{m,\ell}(z,\alpha)$ is a rational function of type $(m,\ell)$.
\item On the interval  $[\alpha^2,1]$, $r_{m,\ell}(z,\alpha)/\sqrt{z}$ takes values in $[1-\varepsilon,1+\varepsilon]$ and achieves its extremal values at exactly $m+\ell+2$ points in $[\alpha^2,1]$ in an alternating fashion.
\item On the interval $[\alpha^2,1]$, $\hat{r}_{m,\ell}(z,\alpha)/\sqrt{z}$ takes values in $[1,(1+\varepsilon)/(1-\varepsilon)]$.
\end{enumerate}
It follows, in particular, that
\[
\max_{z \in [\alpha^2,1]} |r_{m,\ell}(z,\alpha)/\sqrt{z}-1| = \varepsilon_{m,\ell}(\alpha).
\]
Note that the scaled function $\hat{r}_{m,\ell}(z,\alpha)$ is related to the scaled function $\hat{s}_{2\ell+1,2m}(z,\alpha)$ in a simple way:
\begin{equation} \label{shatrhat}
\hat{s}_{2\ell+1,2m}(z,\alpha) =  \frac{z}{\hat{r}_{m,\ell}(z^2,\alpha)}.
\end{equation}

\paragraph{Error estimates}

The errors~(\ref{errsign}) are known to satisfy
\[
E_{p,p}(\alpha) = \frac{2\sqrt{Z_p(\alpha)}}{1+Z_p(\alpha)}
\]
for each $p$, where
\[
Z_p(\alpha) = \inf_{r \in \mathcal{R}_{p,p}} \frac{\sup_{z \in [\alpha,1]}|r(z)|}{\inf_{z \in [-1,-\alpha]}|r(z)|}
\]
is the \emph{Zolotarev number} of the sets $[-1,-\alpha]$ and $[\alpha,1]$~\cite[p. 9]{beckermann2017singular}.  An explicit formula for $Z_p(\alpha)$ is given in~\cite[Theorem 3.1]{beckermann2017singular}.  For our purposes, it is enough to know that $Z_p(\alpha)$ obeys an asymptotically sharp inequality~\cite[Corollary 3.2]{beckermann2017singular}
\[
Z_p(\alpha) \le 4 \rho(\alpha)^{-2p},
\]
where $\rho(\alpha)$ is given by~(\ref{rho}).
This, together with the fact that $E_{p,p}=E_{2\floor{(p-1)/2}+1,2\floor{p/2}}$ for every $p$~\cite[p. 22]{beckermann2017singular}, shows that
\begin{align}
\varepsilon_{m,m-1}(\alpha) &= E_{2m-1,2m}(\alpha) = E_{2m,2m}(\alpha) \le 2\sqrt{Z_{2m}(\alpha)} \le 4 \rho(\alpha)^{-2m}, \label{epsilonmm-1} \\
\varepsilon_{m,m}(\alpha) &= E_{2m+1,2m}(\alpha) = E_{2m+1,2m+1}(\alpha) \le 2\sqrt{Z_{2m+1}(\alpha)} \le 4 \rho(\alpha)^{-(2m+1)}, \label{epsilonmm}
\end{align}
and these bounds are asymptotically sharp.  (The upper bound for $\varepsilon_{m,m-1}(\alpha)$ also appears in~\cite[p. 151, Theorem 5.5]{braess1986nonlinear}.)

\subsection{Proofs}

\paragraph{Proof of Theorem~\ref{thm:rationalzolo}}
To prove Theorem~\ref{thm:rationalzolo}, it will be convenient to introduce the notation 
\[
\hat{s}_{2m} := \hat{s}_{2m-1,2m}, \quad \hat{r}_{2m} := \hat{r}_{m,m-1}, \quad \varepsilon_{2m} := \varepsilon_{m,m-1},
\]
and
\[
\hat{s}_{2m+1} := \hat{s}_{2m+1,2m}, \quad \hat{r}_{2m+1} := \hat{r}_{m,m}, \quad \varepsilon_{2m+1} := \varepsilon_{m,m}.
\]
In this notation, the relation~(\ref{shatrhat}) takes the form
\[
\hat{s}_p(z,\alpha) = \frac{z}{\hat{r}_p(z^2,\alpha)}
\]
for every $p$.  Equivalently,
\[
\hat{r}_p(z,\alpha) = \frac{\sqrt{z}}{\hat{s}_p(\sqrt{z},\alpha)}.
\]
In addition, the relation~(\ref{rhat}) takes the form
\begin{equation} \label{rhat2}
\hat{r}_p(z,\alpha) = \frac{1}{1-\varepsilon_p(\alpha)} r_p(z,\alpha).
\end{equation}

It has been shown in~\cite[Corollary 4]{nakatsukasa2016computing} that if $p$ is odd and a sequence of scalars $\alpha_0,\alpha_1,\dots$ is defined inductively by
\begin{equation} \label{alphaupdate}
\alpha_{k+1} = \hat{s}_{p}(\alpha_k,\alpha_k), \quad \alpha_0 = \alpha,
\end{equation}
then
\begin{equation} \label{shatrecursion}
\hat{s}_{p^{k+1}}(z,\alpha) = \hat{s}_p(\hat{s}_{p^k}(z,\alpha), \alpha_k)
\end{equation}
for every $k \ge 1$.  As remarked in~\cite{nakatsukasa2016computing}, a nearly identical proof shows that this holds also if $p$ is even.

Now suppose that $f_k(z)$ is defined recursively by
\begin{align} 
f_{k+1}(z) &= f_k(z) \hat{r}_p\left( \frac{z}{f_k(z)^2}, \alpha_k \right), & f_0(z) &= 1, \label{rationalzoloAB1} \\
\alpha_{k+1} &= \frac{\alpha_k}{\hat{r}_p(\alpha_k^2,\alpha_k)}, & \alpha_0 &= \alpha. \label{rationalzoloAB2}
\end{align}
We will show that $f_k(z) = \hat{r}_{p^k}(z,\alpha) = \sqrt{z}/\hat{s}_{p^k}(\sqrt{z},\alpha)$ for every $k \ge 1$ by induction.  Note that~(\ref{rationalzoloAB2}) generates the same sequence of scalars as~(\ref{alphaupdate}), and clearly $f_1(z) = \hat{r}_{p}(z,\alpha)$.  If $f_k(z) = \hat{r}_{p^k}(z,\alpha)$ for some $k$, then
\begin{align*}
\hat{r}_{p^{k+1}}(z,\alpha)
&= \frac{\sqrt{z}}{\hat{s}_{p^{k+1}}(\sqrt{z},\alpha)} \\
&= \frac{\sqrt{z}}{\hat{s}_p(\hat{s}_{p^k}(\sqrt{z},\alpha),\alpha_k)} \\
&= \frac{\sqrt{z}}{\hat{s}_{p}(\sqrt{z}/f_k(z),\alpha_k)} \\
&= \frac{\sqrt{z}}{\frac{\sqrt{z}/f_k(z)}{\hat{r}_{p}(z/f_k(z)^2,\alpha_k)}} \\
&= f_k(z) \hat{r}_{p}(z/f_k(z)^2,\alpha_k) \\
&= f_{k+1}(z),
\end{align*}
as desired.  

If $p$ is even, i.e. $p=2m$ for some $m$, then $\hat{r}_p = \hat{r}_{m,m-1}$, and~(\ref{rationalzoloAB1}-\ref{rationalzoloAB2}) is equivalent to~(\ref{rationalzoloA1}-\ref{rationalzoloA2}).  Letting $q = \frac{1}{2}(2m)^k$ so that $2q = (2m)^k = p^k$, we conclude that
\[
f_k(z) = \hat{r}_{2q}(z,\alpha) = \hat{r}_{q,q-1}(z,\alpha), \quad q = \frac{1}{2}(2m)^k.
\]

On the other hand, if $p$ is odd, i.e. $p=2m+1$ for some $m$, then $\hat{r}_p = \hat{r}_{m,m}$, and~(\ref{rationalzoloAB1}-\ref{rationalzoloAB2}) is equivalent to~(\ref{rationalzoloB1}-\ref{rationalzoloB2}).  Letting $q = \frac{1}{2}((2m+1)^k-1)$ so that $2q+1 = (2m+1)^k = p^k$, we conclude that
\[
f_k(z) = \hat{r}_{2q+1}(z,\alpha) = \hat{r}_{q,q}(z,\alpha), \quad q = \frac{1}{2}((2m+1)^k-1).
\]

It remains to prove that $\hat{r}_{p^k}(z,\alpha) = \left(\frac{1+\alpha_k}{2\alpha_k}\right) r_{p^k}(z,\alpha)$ for every $k \ge 1$.  In view of~(\ref{rhat2}), this is equivalent to proving that
\begin{equation} \label{alphaeps}
\alpha_k = \frac{1-\varepsilon_{p^k}(\alpha)}{1+\varepsilon_{p^k}(\alpha)}.
\end{equation}
From~(\ref{equipts}) and~(\ref{shat}), we know that
\[
\hat{s}_{p^k}(\alpha,\alpha) =  \frac{1-\varepsilon_{p^k}(\alpha)}{1+\varepsilon_{p^k}(\alpha)},
\]
so it suffices to show that
\begin{equation} \label{alphakformula}
\alpha_k = \hat{s}_{p^k}(\alpha,\alpha)
\end{equation}
for every $k \ge 1$.
We prove this by induction.  The base case is clear, and if~(\ref{alphakformula}) holds for some $k$, then, upon applying~(\ref{alphaupdate}) and~(\ref{shatrecursion}) with $z=\alpha$, we see that
\begin{align*}
\alpha_{k+1} 
&= \hat{s}_{p}(\alpha_k,\alpha_k) \\
&= \hat{s}_{p}(\hat{s}_{p^k}(\alpha,\alpha), \alpha_k) \\
&= \hat{s}_{p^{k+1}}(\alpha,\alpha).
\end{align*}

\paragraph{Proof of Theorem~\ref{thm:convergence}}

Theorem~\ref{thm:convergence} is a consequence of the following lemma, which we will prove in nearly the same way that Beckermann proves~\cite[Theorem 4]{beckermann2013optimally}.
\begin{lemma} \label{lemma:scalarerr}
Let $m \in \mathbb{N}$, $\ell \in \{m-1,m\}$, $\alpha \in (0,1)$, and $z \in \mathbb{C} \setminus ((-\infty,0] \cup [\alpha^2,1])$.  If $\max\{2|\varphi(z,\alpha)|^{-2(m+\ell+1)},4\rho(\alpha)^{-2(m+\ell+1)}\} < 1$, then
\begin{align*}
|r_{m,\ell}(z,\alpha)/\sqrt{z}-1| &\le \frac{4|\varphi(z,\alpha)|^{-(m+\ell+1)} + 8\rho(\alpha)^{-2(m+\ell+1)}}{\left(1-2|\varphi(z,\alpha)|^{-2(m+\ell+1)}\right)\left(1-4\rho(\alpha)^{-2(m+\ell+1)} \right)} \\&= 4|\varphi(z,\alpha)|^{-(m+\ell+1)} + O\left(|\varphi(z,\alpha)|^{-2(m+\ell+1)}\right).
\end{align*}
where $\varphi(z,\alpha)$ and $\rho(\alpha)$ are given by~(\ref{phi}) and~(\ref{rho}).
\end{lemma}
\paragraph{Remark} When $z \in [\alpha^2,1]$, the slightly sharper bound
\[
|r_{m,\ell}(z,\alpha)/\sqrt{z} - 1| \le 4\rho(\alpha)^{-(m+\ell+1)}
\]
holds in view of~(\ref{epsilonmm-1}-\ref{epsilonmm}).
\begin{proof}

With $Z:=Z_{m+\ell+1}(\alpha)$, let
\[
R(z) = \frac{1 - \left(\frac{1+Z}{1-Z}\right)s_{2\ell+1,2m}(z,\alpha)}{1 + \left(\frac{1+Z}{1-Z}\right)s_{2\ell+1,2m}(z,\alpha)}.
\]
Since $s_{2\ell+1,2m}(z,\alpha)$ takes values in $[1-2\sqrt{Z}/(1+Z),1+2\sqrt{Z}/(1+Z)]$ on the interval $[\alpha,1]$, $R(z)$ takes values in $[-\sqrt{Z},\sqrt{Z}]$ on $[\alpha,1]$.  On the other hand, since $s_{2\ell+1,2m}(z,\alpha)$ is purely imaginary for $z \in i\mathbb{R}$, $|R(z)|=1$ for $z \in i\mathbb{R}$. 

Recall that $\varphi(z,\alpha)$ supplies a conformal map from $\mathbb{C} \setminus ((-\infty,0] \cup [\alpha^2,1])$ to the annulus $\{z \in \mathbb{C} : 1 < |z| < \rho(\alpha)\}$.  Thus, by the maximum principle,
\begin{align*}
\sup_{z \in \mathbb{C} \setminus ((-\infty,0] \cup [\alpha^2,1])} |\varphi(z,\alpha)|^{m+\ell+1} |R(\sqrt{z})|
&=  \sup_{z \in \mathbb{C} \setminus ((-\infty,0] \cup [\alpha^2,1])} |\varphi(z,\alpha)^{m+\ell+1} R(\sqrt{z})| \\
&= \sup_{w \in i\mathbb{R} \cup [\alpha,1]} |\varphi(w^2,\alpha)^{m+\ell+1} R(w)| \\
&\le \max\{1,\rho(\alpha)^{m+\ell+1} \sqrt{Z}\}.
\end{align*}
Since 
\begin{equation} \label{Zbound}
Z = Z_{m+\ell+1}(\alpha) \le 4\rho(\alpha)^{-2(m+\ell+1)},
\end{equation} 
it follows that 
\begin{equation} \label{Rbound}
|R(\sqrt{z})| \le 2|\varphi(z,\alpha)|^{-(m+\ell+1)}
\end{equation}
for every $z \in \mathbb{C} \setminus ((-\infty,0] \cup [\alpha^2,1])$.

Now observe that
\[
r_{m,\ell}(z,\alpha) = \frac{\sqrt{z}}{s_{2\ell+1,2m}(\sqrt{z},\alpha)} = \left( \frac{1+Z}{1-Z} \right)  \left( \frac{1+R(\sqrt{z})}{1-R(\sqrt{z})} \right) \sqrt{z},
\]
so
\begin{align*}
r_{m,\ell}(z,\alpha)/\sqrt{z} - 1 
&= \frac{2(R(\sqrt{z})+Z)}{(1-Z)(1-R(\sqrt{z}))}.
\end{align*}
Invoking~(\ref{Zbound}-\ref{Rbound}) completes the proof.
\end{proof}

\paragraph{Proof of Proposition~\ref{prop:pade}}

It is straightforward to deduce from~\cite[Theorem 5.9]{higham2008functions} the following explicit formula for the type $(m,\ell)$ Pad\'e approximant of $\sqrt{z}$ at $z=1$ for $\ell \in \{m-1,m\}$:
\[
p_{m,\ell}(z) = \sqrt{z} \frac{(1+\sqrt{z})^{m+\ell+1} + (1-\sqrt{z})^{m+\ell+1}}{(1+\sqrt{z})^{m+\ell+1} - (1-\sqrt{z})^{m+\ell+1}}.
\]
It is then easy to check by direct substitution that the roots and poles of $p_{m,\ell}(z)$ are $\left\{ -\tan^2\left(\frac{(2j-1)\pi}{2(m+\ell+1)}\right) \right\}_{j=1}^m$ and $\left\{ -\tan^2\left(\frac{j\pi}{m+\ell+1}\right) \right\}_{j=1}^\ell$, respectively.

On the other hand, the roots and poles of
\[
\hat{r}_{m,\ell}(z,\alpha) = \frac{1+\varepsilon_{m,\ell}(\alpha)}{M(\alpha)} \frac{\prod_{j=1}^{m} (z+c_{2j-1}(\alpha))}{\prod_{j=1}^{\ell} (z+c_{2j}(\alpha))}
\]
are $\{c_{2j-1}(\alpha)\}_{j=1}^m$ and $\{c_{2j}(\alpha)\}_{j=1}^\ell$, respectively, where $c_j(\alpha)$ is given by~(\ref{cj}).  These approach the roots and poles of $p_{m,\ell}(z)$, since the identities $K(0)=\pi/2$, $\mathrm{sn}(z,0)=\sin z$, and $\mathrm{cn}(z,0)=\cos z$ imply that
\begin{align*}
\lim_{\alpha \uparrow 1} c_j(\alpha) &= \tan^2\left(\frac{j\pi}{2(m+\ell+1)}\right).
\end{align*}
The proof is completed by noting that $\hat{r}_{m,\ell}$ is scaled in such a way that $\lim_{\alpha \uparrow 1} \hat{r}_{m,\ell}(1,\alpha)=1=p_{m,\ell}(1)$. 

\paragraph{Remark}  
Alternatively, Proposition~\ref{prop:pade} can be proved by appealing to a general result concerning the convergence of minimax approximants to Pad\'e approximants~\cite{trefethen1985convergence}.  We showed above that $p_{m,\ell}(z)$ is nondegenerate (it has exactly $m$ roots and $\ell$ poles), so Theorem 3b of~\cite{trefethen1985convergence} implies that $r_{m,\ell}(z,\alpha)$ (and hence $\hat{r}_{m,\ell}(z,\alpha)$) converges coefficient-wise to $p_{m,\ell}(z)$ as $\alpha \uparrow 1$.

\paragraph{Proof of Theorem~\ref{thm:stability}}

It is clear from~(\ref{alphaeps}) and~(\ref{epsilonmm-1}-\ref{epsilonmm}) that $\alpha_k \rightarrow 1$ with order of convergence $m+\ell+1$.  Now let $A \in \mathbb{C}^{n \times n}$ have no nonpositive real eigenvalues.  For $\epsilon > 0 $ sufficiently small, the pseudo-spectrum $\Lambda_\epsilon(A)$ is compactly contained in $\mathbb{C} \setminus (\infty,0]$, so $\sup_{z \in \overline{\Lambda_\epsilon(A)}} |\varphi(z,\alpha)| > 1$.  Since $\overline{\Lambda_\epsilon(A)}$ is a spectral set for $A$, we conclude from Theorem~\ref{thm:convergence} that $X_k \rightarrow A^{1/2}$ with order of convergence $m+\ell+1$ in the iteration~(\ref{zolo1}-\ref{zolo2}).  Since $Y_k = X_k^{-1} A$ and $Z_k = X_k^{-1}$ in~(\ref{zolocoupled1}-\ref{zolocoupled3}), it follows that $Y_k \rightarrow A^{1/2}$ and $Z_k \rightarrow A^{-1/2}$ with order of convergence $m+\ell+1$.  
Stability follows easily from the fact that~(\ref{zolocoupled1}-\ref{zolocoupled3}) reduces to the stable Pad\'e iteration~(\ref{padecoupled1}-\ref{padecoupled2}) as $\alpha_k \rightarrow 1$.  (Indeed, in finite-precision arithmetic, there is an integer $K$ such that the computed $\alpha_k$ is rounded to $1$ for every $k \ge K$.)

\section{Practical Considerations} \label{sec:practical}

In this section, we discuss the implementation of the Zolotarev iterations, strategies for terminating the iterations, and computational costs.

\subsection{Implementation}

To implement the Zolotarev iteration~(\ref{zolocoupled1}-\ref{zolocoupled3}), we advocate the use of a partial fraction expansion of $h_{\ell,m}(\cdot,\alpha)$. since it enhances parellelizability and, in our experience, tends to improve stability.   In partial fraction form, $h_{\ell,m}(z,\alpha) = \hat{r}_{m,\ell}(z,\alpha)^{-1}$ for $\ell \in \{m-1,m\}$ is given by
\begin{align}
h_{m-1,m}(z,\alpha) &= \hat{M}(\alpha) \sum_{j=1}^m \frac{a_j(\alpha)}{z+c_{2j-1}(\alpha)}, \label{hm-1m} \\
h_{m,m}(z,\alpha) &= \hat{N}(\alpha) \left( 1 + \sum_{j=1}^m \frac{a_j(\alpha)}{z+c_{2j-1}(\alpha)} \right), \label{hmm}
\end{align}
where
\[
a_j(\alpha) = -\prod_{p=1}^{\ell} (c_{2p}(\alpha) - c_{2j-1}(\alpha)) \bigg/ \prod_{\substack{p=1\\p\neq j}}^{m} (c_{2p-1}(\alpha)-c_{2j-1}(\alpha)),
\]
and $\hat{M}(\alpha)$ and $\hat{N}(\alpha)$ are scalars determined uniquely by the condition that
\[
\min_{z \in [\alpha^2,1]} (h_{\ell,m}(z,\alpha)^{-1}/\sqrt{z} - 1) = 0.
\]
For the reader's benefit, we recall here that
\[
c_j(\alpha) = \alpha^2 \frac{\mathrm{sn}^2\left( \frac{jK(\alpha')}{m+\ell+1}; \alpha' \right)}{\mathrm{cn}^2 \left( \frac{jK(\alpha')}{m+\ell+1}; \alpha' \right)},
\]
where $\mathrm{sn}(\cdot;\alpha)$, $\mathrm{cn}(\cdot;\alpha)$, and $\mathrm{dn}(\cdot;\alpha)$ denote Jacobi's elliptic functions with modulus $\alpha$, $K(\alpha) = \int_0^{\pi/2} (1-\alpha^2 \sin^2\theta)^{-1/2} \, d\theta$ is the complete elliptic integral of the first kind, and $\alpha' = \sqrt{1-\alpha^2}$ is the complementary modulus to $\alpha$.  Note that $c_j(\alpha)$ and $a_j(\alpha)$ depend implicitly on $m$ and $\ell$.  In particular, $a_j(\alpha)$ and $c_{2j-1}(\alpha)$ have different values in~(\ref{hm-1m}) (where $\ell=m-1$) than they do in~(\ref{hmm}) (where $\ell=m$).

Since the locations of the minima of $h_{\ell,m}(z,\alpha)^{-1}/\sqrt{z} = \hat{r}_{m,\ell}(z,\alpha)/\sqrt{z} = \linebreak \hat{s}_{2\ell+1,2m}(\sqrt{z},\alpha)^{-1}$ follow from~(\ref{equipts}), one can obtain explicit expressions for $\hat{M}(\alpha)$ and $\hat{N}(\alpha)$:
\begin{align*}
\hat{M}(\alpha) &= \left( \sqrt{\zeta} \sum_{j=1}^m \frac{a_j(\alpha)}{\zeta+c_{2j-1}(\alpha)} \right)^{-1}, \quad \zeta = \alpha^2 \bigg/ \mathrm{dn}^2\left( \frac{K(\alpha')}{2m}; \alpha' \right), \\
\hat{N}(\alpha) &= \left( 1 + \sum_{j=1}^m \frac{a_j(\alpha)}{1+c_{2j-1}(\alpha)} \right)^{-1}.
\end{align*}

Note that accurate evaluation of $K(\alpha')$, $\mathrm{sn}(\cdot;\alpha')$, $\mathrm{cn}(\cdot;\alpha')$, and $\mathrm{dn}(\cdot;\alpha')$ in floating point arithmetic is a delicate task when $\alpha' \approx 1 \iff \alpha \approx 0$~\cite[Section 4.3]{nakatsukasa2016computing}.  Rather than using the built-in MATLAB functions \verb$ellipj$ and \verb$ellipke$ to evaluate these elliptic functions, we recommend using the code described in~\cite[Section 4.3]{nakatsukasa2016computing}, which is tailored for our application.

Written in full, the Zolotarev iteration~(\ref{zolocoupled1}-\ref{zolocoupled3}) of type $(m,m-1)$\footnote{Although $h_{m-1,m}(z,\alpha)$ is a rational function of type $(m-1,m)$, we continue to refer to this iteration as the type $(m,m-1)$ Zolotarev iteration since $\hat{r}_{m,m-1}(z,\alpha) = h_{m-1,m}(z,\alpha)^{-1}$ is of type $(m,m-1)$.} reads
\begin{align}
Y_{k+1} &= \hat{M}(\alpha_k) \sum_{j=1}^m a_j(\alpha_k) Y_k (Z_k Y_k+c_{2j-1}(\alpha_k)I)^{-1}, & Y_0 &= A, \label{zoloA1full} \\
Z_{k+1} &= \hat{M}(\alpha_k) \sum_{j=1}^m a_j(\alpha_k) (Z_k Y_k+c_{2j-1}(\alpha_k)I)^{-1} Z_k, & Z_0 &= I, \label{zoloA2full} \\
\alpha_{k+1} &= \alpha_k h_{m-1,m}(\alpha_k^2,\alpha_k), & \alpha_0 &= \alpha, \label{zoloA3full}
\end{align}
and the Zolotarev iteration of type $(m,m)$ reads
\begin{align}
Y_{k+1} &= \hat{N}(\alpha_k) \left( Y_k + \sum_{j=1}^m a_j(\alpha_k) Y_k (Z_k Y_k+c_{2j-1}(\alpha_k)I)^{-1} \right), & Y_0 &= A, \label{zoloB1full} \\
Z_{k+1} &= \hat{N}(\alpha_k) \left(Z_k + \sum_{j=1}^m a_j(\alpha_k) (Z_k Y_k+c_{2j-1}(\alpha_k)I)^{-1} Z_k \right), & Z_0 &= I, \label{zoloB2full} \\
\alpha_{k+1} &= \alpha_k h_{m,m}(\alpha_k^2,\alpha_k), & \alpha_0 &= \alpha. \label{zoloB3full}
\end{align}
As alluded to earlier, a suitable choice for $\alpha$ is $\alpha = \sqrt{|\lambda_{\mathrm{min}}(A) / \lambda_{\mathrm{max}}(A)|}$ (or an estimate thereof), and it important to scale $A$ so that its spectral radius is $1$ (or approximately 1).

\subsection{Floating Point Operations}

The computational costs of the Zolotarev iterations depend on the precise manner in which they are implemented.  One option is compute $Z_k Y_k$ (1 matrix multiplication), obtain $h_{\ell,m}(Z_k Y_k,\alpha_k)$ by computing $(Z_k Y_k + c_{2j-1}(\alpha_k)I)^{-1}$ for each $j$ ($m$ matrix inversions), and multiply $Y_k$ and $Z_k$ by $h_{\ell,m}(Z_k Y_k,\alpha_k)$ (2 matrix multiplications). 
An alternative that is better suited for parallel computations is to compute $Z_k Y_k$ (1 matrix multiplication), compute the LU factorization $L_j U_j = Z_k Y_k + c_{2j-1}(\alpha_k)I$ for each $j$ ($m$ LU factorizations), and perform $m$ ``right divisions by a factored matrix'' $Y_k (L_j U_j)^{-1}$ and $m$ ``left divisions by a factored matrix'' $(L_j U_j)^{-1} Z_k$ via forward and back substitution.  
In parallel, all $m$ LU factorizations can be performed simultaneously, and all $2m$ divisions by factored matrices can be performed simultaneously, so that the effective cost per iteration is $\frac{14}{3}n^3$ flops.  In the first iteration, the cost reduces to $\frac{8}{3}n^3$ flops since $Z_0=I$.  The total effective cost for $k$ iterations is $(\frac{8}{3}+\frac{14}{3}(k-1))n^3$ flops, which is less than the (serial) cost of a direct method, $28\frac{1}{3}n^3$ flops~\cite[p. 136]{higham2008functions}, whenever $k \le 6$.

Yet another alternative is to write~(\ref{zoloA1full}-\ref{zoloA2full}) in the form
\begin{align}
Y_{k+1} &= \hat{M}(\alpha_k) \left( \sum_{j=1}^m a_j(\alpha_k) Y_k (Y_k+c_{2j-1}(\alpha_k)Z_k^{-1})^{-1} \right) Z_k^{-1}, & Y_0 &= A, \label{zoloA1alt} \\
Z_{k+1} &= \hat{M}(\alpha_k) \sum_{j=1}^m a_j(\alpha_k) (Y_k+c_{2j-1}(\alpha_k)Z_k^{-1})^{-1}, & Z_0 &= I, \label{zoloA2alt}
\end{align}
and similarly for~(\ref{zoloB1full}-\ref{zoloB2full}).
Interestingly, this form of the iteration has exhibited the best accuracy in our numerical experiments, for reasons that are not well understood.  It can be parallelized by performing the $m$ right divisions $Y_k (Y_k+c_{2j-1}(\alpha_k)Z_k^{-1})^{-1}$ and $m$ inversions $(Y_k+c_{2j-1}(\alpha_k)Z_k^{-1})^{-1}$ simultaneously, recycling LU factorizations in the obvious way.  Moreover, the final multiplication by $Z_k^{-1}$ in~(\ref{zoloA1alt}) can be performed in parallel with the inversion of $Z_{k+1}$.  The effective cost in such a parallel implementation is $\frac{14}{3}kn^3$ flops.

\subsection{Termination Criteria}

We now consider the question of how to terminate the iterations.  Define $\widetilde{X}_k = 2\alpha_k X_k / (1+\alpha_k)$, $\widetilde{Y}_k = (1+\alpha_k) Y_k / (2\alpha_k)$, and $\widetilde{Z}_k = (1+\alpha_k) Z_k / (2\alpha_k)$.  Since $\widetilde{X}_k$, $\widetilde{Y}_k$, $\widetilde{Z}_k$, and $A$ commute with one another, and since $\widetilde{Y}_k = \widetilde{X}_k^{-1} A$ and $\widetilde{Z}_k = \widetilde{X}_k^{-1} = \widetilde{Y}_k A^{-1}$, it is easy to verify that
\[
(\widetilde{Y}_k A^{-1/2}-I) + (\widetilde{Z}_k A^{1/2}-I) = (\widetilde{Z}_k \widetilde{Y}_k - I) - (\widetilde{Z}_k A^{1/2}-I) (\widetilde{Y}_k A^{-1/2}-I)
\]
and
\[
\widetilde{Y}_k A^{-1/2}-I = (I - \widetilde{X}_k A^{-1/2}) + (\widetilde{Y}_k A^{-1/2}-I) (I - \widetilde{X}_k A^{-1/2}).
\]
By dropping second order terms, we see that near convergence,
\begin{equation} \label{err_relations}
I - \widetilde{X}_k A^{-1/2}.\approx \widetilde{Y}_k A^{-1/2}-I = \widetilde{Z}_k A^{1/2}-I \approx \frac{1}{2} (\widetilde{Z}_k \widetilde{Y}_k - I).
\end{equation}
The relative errors $\frac{\|\widetilde{Y}_k-A^{1/2}\|}{\|A^{1/2}\|} \le \|\widetilde{Y}_k A^{-1/2}-I\|$ and $\frac{\|\widetilde{Z}_k-A^{-1/2}\|}{\|A^{-1/2}\|} \le \|\widetilde{Z}_k A^{1/2}-I\|$ will therefore be (approximately) smaller than a tolerance $\delta>0$ so long as
\begin{equation} \label{ZY-I}
\|\widetilde{Z}_k \widetilde{Y}_k - I\| \le 2\delta.
\end{equation}

While theoretically appealing, the criterion~(\ref{ZY-I}) is not ideal for computations for two reasons.  It costs an extra matrix multiplication in the last iteration, and, more importantly,~(\ref{ZY-I}) may never be satisfied in floating point arithmetic.  A cheaper, more robust option is to approximate $\|\widetilde{Z}_k \widetilde{Y}_k - I\|$ based on the value of $\|\widetilde{Z}_{k-1} \widetilde{Y}_{k-1} - I\|$ as follows. In view of~(\ref{err_relations}) and Theorem~\ref{thm:convergence}, we have
\[
\|\widetilde{Z}_k \widetilde{Y}_k - I\|_2 \lessapprox 8 L \gamma^{-(m+\ell+1)^k}
\]
for some constants $L \ge 1$ and $\gamma > 1$.  Denoting $e_k := 8 L \gamma^{-(m+\ell+1)^k}$, we have
\[
e_k \le 2\delta \iff e_{k-1} \le 8L \left( \frac{\delta}{4L} \right)^{1/(m+\ell+1)}.
\]
This suggests that we terminate the iteration and accept $\widetilde{Y}_k$ and $\widetilde{Z}_k$ as soon as
\[
\|\widetilde{Z}_{k-1} \widetilde{Y}_{k-1} - I\|_2 \le 8L \left( \frac{\delta}{4L} \right)^{1/(m+\ell+1)}.
\]
In practice, $L$ is not known, and it may be preferable to use a different norm, so we advocate terminating when 
\begin{equation} \label{termination1}
\|\widetilde{Z}_{k-1} \widetilde{Y}_{k-1} - I\| \le 8 \left( \frac{\delta}{4} \right)^{1/(m+\ell+1)},
\end{equation}
where $\delta$ is a relative error tolerance with respect to a desired norm $\|\cdot\|$.  Note that this test comes at no additional cost if the product $\widetilde{Z}_{k-1} \widetilde{Y}_{k-1}$ was computed at iteration $k-1$.  
If $\widetilde{Z}_{k-1}^{-1}$ is known but $\widetilde{Z}_{k-1} \widetilde{Y}_{k-1}$ is not (as is the case when~(\ref{zoloA1alt}-\ref{zoloA2alt}) is used), then we have found the following criterion, which is inspired by~\cite[Equation 6.31]{higham2008functions}, to be an effective alternative:
\begin{equation} \label{termination2}
\|\widetilde{Y}_k-\widetilde{Y}_{k-1}\| \le \left( \delta \frac{\|\widetilde{Y}_k\|}{\|A^{-1}\|\|\widetilde{Z}_{k-1}^{-1}\|} \right)^{1/(m+\ell+1)}.
\end{equation}
In either case, we recommend also terminating the iteration if the relative change in $\widetilde{Y}_k$ is small but fails to decrease significantly, e.g.,
\begin{equation} \label{reldiff}
\frac{1}{2} \frac{\|\widetilde{Y}_{k-1}-\widetilde{Y}_{k-2}\|}{\|\widetilde{Y}_{k-1}\|} \le \frac{\|\widetilde{Y}_k-\widetilde{Y}_{k-1}\|}{\|\widetilde{Y}_k\|} \le 10^{-2}.
\end{equation}

\section{Numerical Examples} \label{sec:numerical}

In this section, we study the performance of the Zolotarev iterations with numerical experiments. 

\subsection{Scalar Iteration} \label{sec:scalar}

\begin{figure}
\includegraphics[page=1,trim=0.9in 2.5in 0.9in 2.5in,clip,scale=0.75]{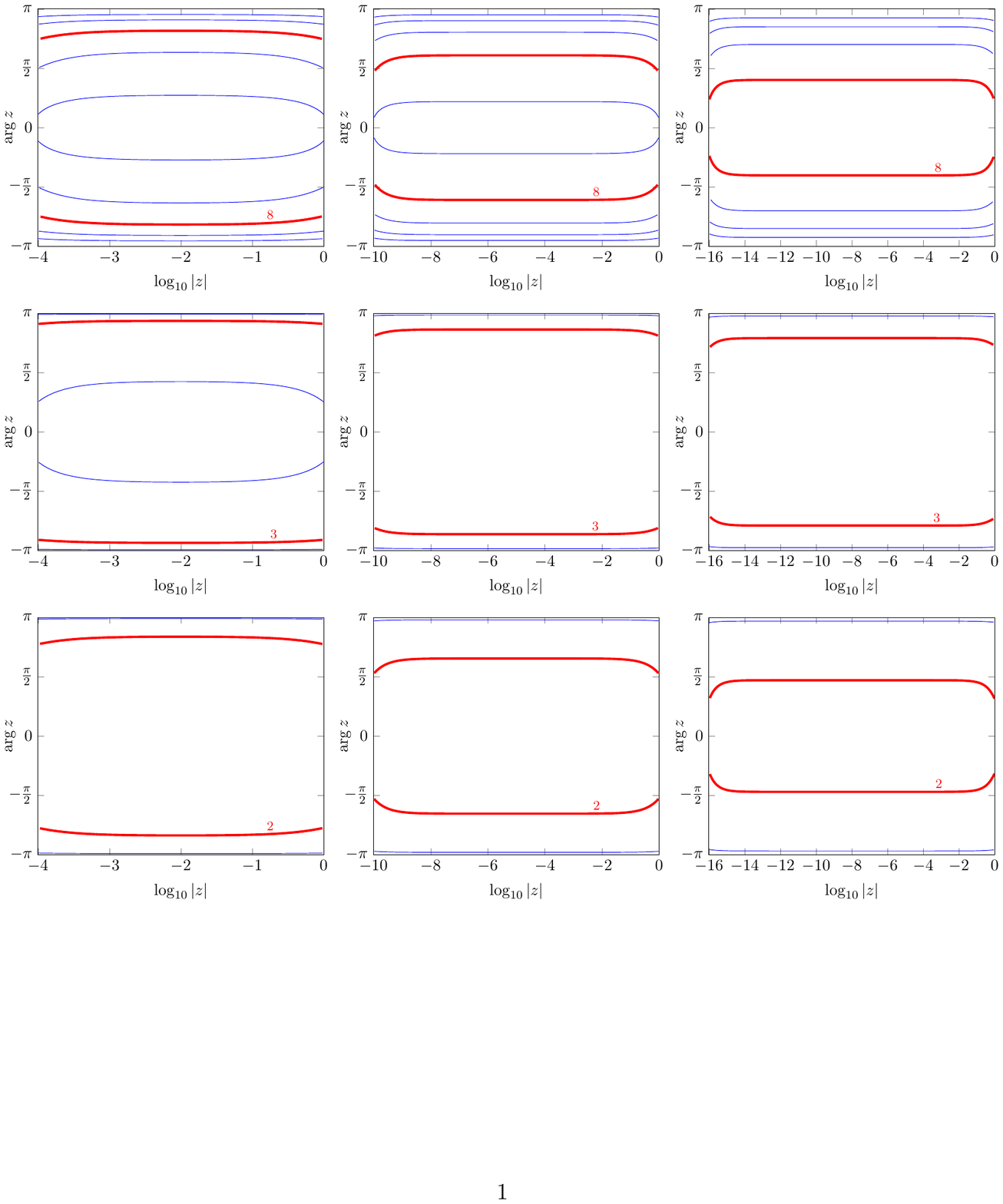}
\caption{Integer level sets of $\kappa(z,\alpha)$ for $(m,\ell) = (1,0)$ (row 1), $(m,\ell) = (4,4)$ (row 2), $(m,\ell) = (8,8)$ (row 3), $\alpha = 10^{-2}$ (column 1), $\alpha = 10^{-5}$ (column 2), and $\alpha = 10^{-8}$ (column 3).  To help compare level sets within each row, we have arbitrarily selected a single level set to label and highlight in bold red.  Each unlabelled level set's value differs from that of its nearest inner neighbor by $+1$.} 
\label{fig:contours}
\end{figure}

\begin{figure}
\includegraphics[page=2,trim=0.9in 2.5in 0.9in 2.5in,clip,scale=0.75]{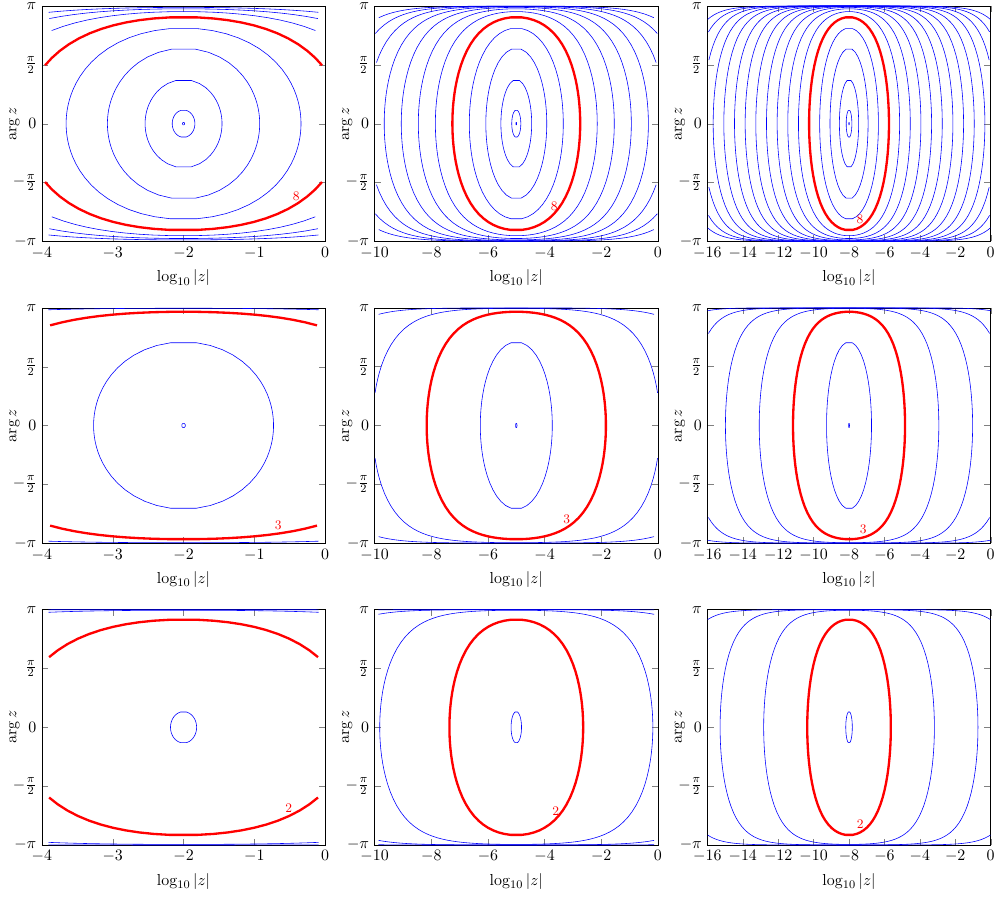}
\caption{Integer level sets of $\kappa(z/\alpha,1)$ for $(m,\ell) = (1,0)$ (row 1), $(m,\ell) = (4,4)$ (row 2), $(m,\ell) = (8,8)$ (row 3), $\alpha = 10^{-2}$ (column 1), $\alpha = 10^{-5}$ (column 2), and $\alpha = 10^{-8}$ (column 3).  To help compare level sets within each row, we have arbitrarily selected a single level set to label and highlight in bold red.  Each unlabelled level set's value differs from that of its nearest inner neighbor by $+1$.} 
\label{fig:contourspade}
\end{figure}

To gain some intuition behind the behavior of the Zolotarev iteration for the matrix square root, we begin by investigating the behavior of its scalar counterpart.

Lemma~\ref{lemma:scalarerr} shows that if $f_k(z)$ and $\alpha_k$ are defined as in Theorem~\ref{thm:rationalzolo}, then
\begin{equation}
\left| \left(\frac{2\alpha_k}{1+\alpha_k}\right) f_k(z)/\sqrt{z}-1 \right| \le 4|\varphi(z,\alpha)|^{-(m+\ell+1)^k} + O\left( |\varphi(z,\alpha)|^{-2(m+\ell+1)^k} \right),
\end{equation}
where $\ell=m-1$ if~(\ref{rationalzoloA1}-\ref{rationalzoloA2}) is used, and $\ell=m$ if~(\ref{rationalzoloB1}-\ref{rationalzoloB2}) is used.  Thus, for a given $z \in \mathbb{C} \setminus (\infty,0]$ and a given relative tolerance $\delta > 0$, we can estimate the smallest $k$ for which $|2\alpha_k f_k(z)/((1+\alpha_k)\sqrt{z})-1| \le \delta$: we have $k \approx \lceil \kappa(z,\alpha) \rceil$ with 
\[
\kappa(z,\alpha) = \frac{\log\log(4/\delta) - \log\log|\varphi(z,\alpha)|}{\log(m+\ell+1)}.
\]

Fig.~\ref{fig:contours} plots the integer level sets of $\kappa(z,\alpha)$ for $(m,\ell) \in \{(1,0),(4,4),(8,8)\}$, $\delta = 10^{-16}$, and $\alpha \in \{10^{-2},10^{-5},10^{-8}\}$ in the slit annulus $\mathcal{A} = \{z \in \mathbb{C} \mid \alpha^2 \le |z| \le 1, \, -\pi < \arg z < \pi\}$.  To improve the clarity of the plots, we have plotted the level sets in the $(\log_{10}|z|,\arg z)$ coordinate plane rather than the usual $(\Re z, \Im z)$ coordinate plane.  The level sets have the following interpretation: 
If $z_0 \in \mathbb{C}$ lies within the region enclosed by the level set $\kappa(z,\alpha)=c \in \mathbb{N}$, then the sequence $\{2\alpha_k f_k(z_0)/(1+\alpha_k)\}_{k=0}^\infty$ generated by the type $(m,\ell)$ Zolotarev iteration from Theorem~\ref{thm:rationalzolo} converges to $\sqrt{z_0}$ in at most $c$ iterations with a relative tolerance of $\approx 10^{-16}$.

Observe that when $(m,\ell)=(8,8)$ and $z_0$ lies in the right half-annulus $\{z : \Re z \ge 0, \, \alpha^2 \le |z| \le 1 \}$ (which corresponds to the horizontal strip $\{z : 2\log_{10} \alpha \le \log_{10} |z| \le 0, -\pi/2 < \arg z < \pi/2 \}$ in Fig.~\ref{fig:contours}), convergence of the scalar iteration is achieved in just 2 iterations whenever $\alpha \ge 10^{-5}$.  For nearly all other $z_0 \in \mathcal{A}$, 3 iterations suffice.

\paragraph{Comparison with Pad\'e iterations}

For comparison, Fig.~\ref{fig:contourspade} plots the integer level sets of $\kappa(z/\alpha,1)$ for the same values of $(m,\ell)$, $\delta$, and $\alpha$ as above.  In view of Proposition~\ref{prop:pade}, these level sets dictate the convergence of the type $(m,\ell)$ Pad\'e iteration with the initial iterate $z_0$ scaled by $1/\alpha$.  For $\alpha = 10^{-2}$ (the leftmost column), the behavior of the Pad\'e iteration is not significantly different from the behavior of the Zolotarev iteration.  However, as $\alpha$ decreases, a clear pattern emerges.  The level sets $\kappa(z/\alpha,1)=c$ do not begin to enclose scalars $z$ with extreme magnitudes ($|z| \approx \alpha^2$ and $|z| \approx 1$) until $c$ is relatively large.  For example, when $\alpha=10^{-8}$ and $(m,\ell)=(8,8)$, the smallest integer $c$ for which the level set $\kappa(z/\alpha,1)=c$ encloses both $z = \alpha^2$ and $z=1$ is $c=5$ (see the lower right plot of~Fig.~\ref{fig:contourspade}).  In contrast, for the Zolotarev iteration with the same $(m,\ell)$ and $\alpha$, the smallest integer $c$ for which $\kappa(z,\alpha)=c$ encloses both $z=\alpha^2$ and $z=1$ is $c=2$ (see the lower right plot of Fig.~\ref{fig:contours}).  The situation is similar when $|z| \in \{\alpha^2,1\}$ and $z$ has nonzero imaginary part.

\paragraph{Implications} 

The preceding observations have important implications for computing the square root of a matrix $A \in \mathbb{C}^{n \times n}$ with no nonpositive real eigenvalues.  Without loss of generality, we may assume that $A$ has been scaled in such a way that its spectrum $\Lambda(A)$ is contained the slit annulus  $\{z \in \mathbb{C} \mid \alpha^2 \le |z| \le 1, \, -\pi < \arg z < \pi\}$ for some $\alpha \in (0,1)$.  Then, if $A$ is normal, the number of iterations needed for the Zolotarev iteration of type $(m,\ell)$ to converge to $A^{1/2}$ (i.e. $\|2\alpha_k X_k A^{-1/2}/(1+\alpha_k)-I\| \lessapprox 10^{-16}$) in exact arithmetic is given by the smallest integer $c$ for which the level set $\kappa(z,\alpha)=c$ encloses $\Lambda(A)$.  For the Pad\'e iteration (with $A$ rescaled by $1/\alpha$) the same statement holds with $\kappa(z,\alpha)$ replaced by $\kappa(z/\alpha,1)$.  

We conclude from the preceding discussion that the Zolotarev iterations are often preferable when $A$ has eigenvalues with widely varying magnitudes (assuming $A$ is normal).  For instance, if $|\lambda_{\mathrm{max}}(A)| / |\lambda_{\mathrm{min}}(A)| =\alpha^{-2} \le 10^{10}$ and the spectrum of $A$ lies in the right half plane, then the Zolotarev iteration of type $(8,8)$ converges in at most 2 iterations, whereas the Pad\'e iteration of type $(8,8)$ converges in at most 4 (see row 3, columns 1-2 of Figs.~\ref{fig:contours}-\ref{fig:contourspade}).
When considering non-normal $A$ and/or the effects of round-off errors, the situation is of course more difficult to analyze, but we address this situation with numerical experiments in Section~\ref{sec:numerical_matrix}.

Note that in the Pad\'e iteration~(\ref{pade}), it is common to scale not only the initial iterate $X_0$, but also subsequent iterates $X_k$, by $\mu_k = |\det (X_k)^{-1/n}|$.
(More precisely, this is accomplished in a mathematically equivalent, numerically stabler way by scaling $Y_k$ and $Z_k$ by $\mu_k^{-1} = |(\det Y_k \det Z_k)^{-1/(2n)}|$ in~(\ref{padecoupled1}-\ref{padecoupled2})~\cite[Equation (3.2)]{higham1997stable}).  These scalars will of course depend on the distribution of the eigenvalues of $A$, but in the case in which $m=\ell$ and $A$ has real eigenvalues with logarithms uniformly distributed in $[2\log_{10}\alpha,0]$, one finds that $\mu_0 = 1/\alpha$ and $\mu_k=1$ for $k \ge 1$, showing that Fig.~\ref{fig:contourspade} is a fair representation of the behavior of the scaled Pad\'e iteration.

\subsection{Matrix Iteration} \label{sec:numerical_matrix}

\begin{table}
\centering
\begin{filecontents*}{matstats.dat}
   1.3934901e+00   1.0979654e+00   1.3989099e+08   2.8108626e+00
   4.0179233e+01   8.3226724e+04   5.6871145e+07   5.2026489e+06
   8.0407976e+01   2.0426296e+05   3.0510106e+10   3.9148769e+06
\end{filecontents*}
\pgfplotstabletypeset[
every head row/.style={before row=\toprule,after row=\midrule},
every last row/.style={after row=\bottomrule},
columns={leftcol,0,1,2,3},
create on use/leftcol/.style={create col/set list={$\alpha_\infty(A^{1/2})$,$\kappa_{\mathrm{sqrt}}(A)$,$\kappa_2(A^{1/2})$}},
columns/leftcol/.style={string type,column type/.add={|}{|},column name={}},
columns/0/.style={sci,sci e,sci zerofill,precision=1,verbatim,column type/.add={}{|},column name={$A_1$}},
columns/1/.style={sci,sci e,sci zerofill,precision=1,verbatim,column type/.add={}{|},column name={$A_2$}},
columns/2/.style={sci,sci e,sci zerofill,precision=1,verbatim,column type/.add={}{|},column name={$A_3$}},
columns/3/.style={sci,sci e,sci zerofill,precision=1,verbatim,column type/.add={}{|},column name={$A_4$}}
]
{matstats.dat}
\vspace{0.05in}
\caption{Properties of the matrices $A_1$, $A_2$, $A_3$ and $A_4$.}
\label{tab:matstats}
\end{table}

\begin{table}[t]
\centering
\maketable{tab1.dat} 
\quad
\maketable{tab2.dat}\\ \vspace{0.1in}
\maketable{tab3.dat} 
\hspace{0.05in}\quad
\maketable{tab4.dat}\hspace{-0.15in}
\vspace{0.05in}
\caption{Numerical results for $A_1$ (upper left), $A_2$ (upper right), $A_3$ (lower left), and $A_4$ (lower right).  Each table shows the number of iterations $k$, relative error $\|\hat{X}-A^{1/2}\|_{\infty} / \|A^{1/2}\|_{\infty}$, and relative residual $\|\hat{X}^2-A\|_{\infty} / \|A\|_{\infty}$ in the computed square root $\hat{X}$ of $A$.} 
\label{tab:4tests}
\end{table}

In what follows, we compare the Zolotarev iterations of type $(m,\ell)$ (hereafter referred to as Z-$(m,\ell)$) with the following other methods: the Denman-Beavers iteration (DB)~\cite[Equation (6.28)]{higham2008functions} (see also~\cite{denman1976matrix}), the product form of the Denman-Beavers iteration (DBp)~\cite[Equation (6.29)]{higham2008functions}, the incremental Newton iteration (IN)~\cite[Equation (6.30)]{higham2008functions}\footnote{Note that Equation (6.30) in~\cite{higham2008functions} contains a typo in the last line: $E_{k+1} = -\frac{1}{2}E_k X_{k+1}^{-1} E_k$ should read $E_{k+1} = -\frac{1}{2}\widetilde{E}_k X_{k+1}^{-1} \widetilde{E}_k$.} (see also~\cite{meini2004matrix,iannazzo2003note}), the principal Pad\'e iterations of type $(m,\ell)$ (P-$(m,\ell)$)~\cite[Equation (6.34)]{higham2008functions} (see also~\cite{higham1997stable,higham2005functions}), and the MATLAB function \verb$sqrtm$.  In the Pad\'e and Zolotarev iterations, we focus on the iterations of type $(1,0)$, $(4,4)$, and $(8,8)$ for simplicity.

In all of the iterations (except the Zolotarev iterations), we use determinantal scaling (as described in~\cite[Section 6.5]{higham2008functions} and~\cite[Equation (3.2)]{higham1997stable}) until the $\infty$-norm relative change in $X_k$ falls below $10^{-2}$.  In the Zolotarev iterations, we use $\alpha = \sqrt{|\lambda_{\mathrm{min}}(A)/\lambda_{\mathrm{max}}(A)|}$, and we scale $A$ so that its spectral radius is 1.  In the Zolotarev and Pad\'e iterations, we use the formulation~(\ref{zoloA1alt}-\ref{zoloA2alt}) and its type-$(m,m)$ counterpart, and we terminate the iterations when either~(\ref{termination2}) or~(\ref{reldiff}) is satisfied in the $\infty$-norm with $\delta = u\sqrt{n}$, where $u = 2^{-53}$ is the unit round-off.  To terminate the DB and IN iterations, we use the following termination criterion~\cite[p. 148]{higham2008functions}: 
$\|X_k-X_{k-1}\|_\infty \le (\delta\|X_k\|_\infty/\|X_{k-1}^{-1}\|_\infty)^{1/2}$ or $\frac{1}{2}\|X_{k-1}-X_{k-2}\|_\infty/\|X_{k-1}\|_\infty \le \|X_k-X_{k-1}\|_\infty/\|X_k\|_\infty \le 10^{-2}$.
To terminate the DBp iteration, we replace the first condition by $\|M_k-I\|_\infty \le \delta$, where $M_k$ is the ``product'' matrix in~\cite[Equation (6.29)]{higham2008functions}.  We impose a maximum of 20 iterations for each method.

\paragraph{Four test matrices in detail} We first consider 4 test matrices studied previously in~\cite[Section 6.6]{higham2008functions}:
\begin{enumerate}
\item $A_1 = I + wv^*$, where $w = \left( 1^2 \; 2^2 \; \dots \; n^2 \right)^*$ and $v = \left( 0^2 \; 1^2 \; 2^2 \; \dots \; (n-1)^2 \right)^*$.
\item $A_2 =$ \verb$gallery('moler',16)$.
\item $A_3 =$ \verb$Q*rschur(8,2e2)*Q'$, where \verb$Q=gallery('orthog',8)$ and \verb$rschur$ is a function from the Matrix Computation Toolbox~\cite{Higham:MCT}.
\item $A_4 =$ \verb$gallery('chebvand',16)$.
\end{enumerate}
Table~\ref{tab:matstats} lists some basic information about these matrices, including:
\begin{itemize}
\item The condition number of the $\infty$-norm relative residual of $A^{1/2}$~\cite[Equation (6.4)]{higham2008functions}:
\[
\alpha_\infty(A^{1/2}) = \frac{\|A^{1/2}\|_\infty^2}{\|A\|_\infty}.
\]
\item The Frobenius-norm relative condition number of the matrix square root at $A$~\cite[Equation (6.2)]{higham2008functions}:
\[
\kappa_{\mathrm{sqrt}}(A) = \frac{\|(I \otimes A^{1/2}) + (A^{1/2} \otimes I)\|_2 \|A\|_F}{\|A^{1/2}\|_F}.
\]
\item The 2-norm condition number of $A^{1/2}$:
\[
\kappa_2(A^{1/2}) = \|A^{1/2}\|_2 \|A^{-1/2}\|_2.
\]
\end{itemize}

Table~\ref{tab:4tests} reports the number of iterations $k$, relative error $\|\hat{X}-A^{1/2}\|_{\infty} / \|A^{1/2}\|_{\infty}$, and relative residual $\|\hat{X}^2-A\|_{\infty} / \|A\|_{\infty}$ in the computed square root $\hat{X}$ of $A$ for each method.  (We computed the ``exact'' $A^{1/2}$ using variable precision arithmetic in MATLAB: \verb$vpa(A,100)^(1/2)$.)  In these tests, the Zolotarev and Pad\'e iterations of a given type $(m,\ell)$ tended to produce comparable errors and residuals, but the Zolotarev iterations almost always took fewer iterations to do so.  With the exception of $A_3$, the Zolotarev, Pad\'e, and incremental Newton iterations achieved forward errors less than or comparable to the MATLAB function \verb$sqrtm$.  On $A_3$, \verb$sqrtm$ performed best, but it is interesting to note that the type $(8,8)$ Zolotarev iteration produced the smallest forward error and smallest residual among the iterative methods.

\begin{figure}
\centering
\pgfplotstableread{err.dat}{\tests}
\begin{tikzpicture} 
\begin{semilogyaxis}[width=0.8\textwidth,
minor tick num=1,
ylabel=Relative Error]
\addplot [black,thick] table [x=0, y=1] {\tests};
\addplot+[only marks,mark=*,black,mark options={fill=black}] table [x=0, y=2] {\tests};
\addplot+[only marks,mark=*,black,mark options={fill=white}] table [x=0, y=3] {\tests};
\addplot+[only marks,mark=triangle,blue,mark options={scale=1.5}] table [x=0, y=4] {\tests};
\addplot+[only marks,mark=diamond,green,mark options={scale=1.5}] table [x=0, y=5] {\tests};
\addplot+[only marks,mark=square,red] table [x=0, y=6] {\tests};
\addplot+[only marks,mark=x,olive] table [x=0, y=7] {\tests};
\legend{$u\kappa_{\mathrm{sqrt}}(A)$\\Best iterative\\Worst iterative\\Z-$(1,0)$\\Z-$(4,4)$\\Z-$(8,8)$\\sqrtm\\};
\end{semilogyaxis}
\end{tikzpicture}
\caption{Relative errors committed by each method on 44 tests, ordered by decreasing condition number $\kappa_{\mathrm{sqrt}}(A)$.}
\label{fig:err}
\end{figure}
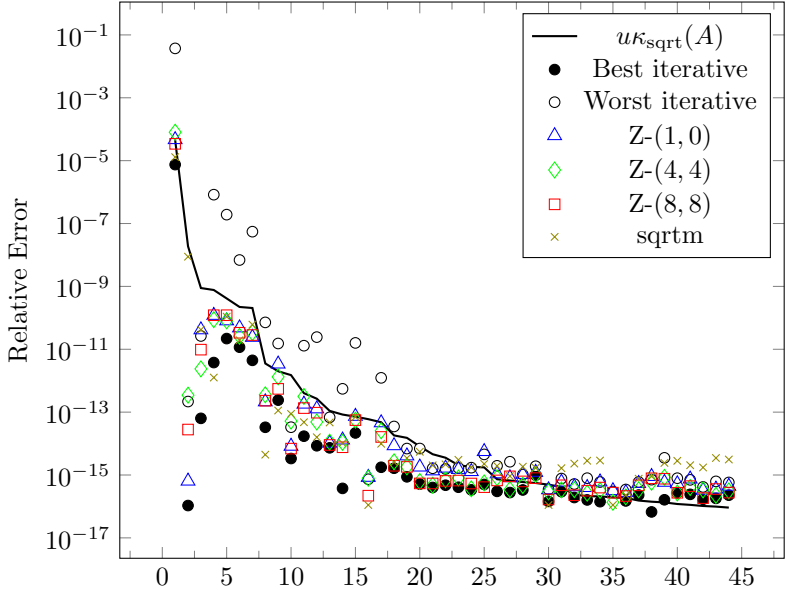

\begin{table}
\centering
\begin{filecontents*}{iter.dat}
   7.3863636e+00   2.1373305e+00   3.0000000e+00   1.2000000e+01
   7.2500000e+00   2.1900940e+00   3.0000000e+00   1.2000000e+01
   7.6818182e+00   2.7851160e+00   3.0000000e+00   2.0000000e+01
   7.6818182e+00   2.3205124e+00   4.0000000e+00   1.3000000e+01
   3.2954545e+00   1.1325942e+00   2.0000000e+00   6.0000000e+00
   2.8181818e+00   9.9470055e-01   2.0000000e+00   5.0000000e+00
   7.5909091e+00   2.0609118e+00   5.0000000e+00   1.2000000e+01
   2.8181818e+00   7.2409454e-01   2.0000000e+00   4.0000000e+00
   2.3636364e+00   4.8660710e-01   2.0000000e+00   3.0000000e+00
   0.0000000e+00   0.0000000e+00   0.0000000e+00   0.0000000e+00
\end{filecontents*}
\pgfplotstabletypeset[
every head row/.style={before row=\toprule,after row=\midrule},
every last row/.style={after row=\bottomrule},
columns={leftcol,0,1,2,3},
create on use/leftcol/.style={create col/set list={DB,DBp,IN,{P-$(1,0)$},{P-$(4,4)$},{P-$(8,8)$},{Z-$(1,0)$},{Z-$(4,4)$},{Z-$(8,8)$},sqrtm}},
columns/leftcol/.style={string type,column type/.add={|}{|},column name={Method}},
columns/0/.style={fixed,precision=1,verbatim,column type/.add={}{|},column name={Mean}},
columns/1/.style={fixed,precision=1,verbatim,column type/.add={}{|},column name={STD}},
columns/2/.style={fixed,precision=1,verbatim,column type/.add={}{|},column name={Min}},
columns/3/.style={fixed,verbatim,column type/.add={}{|},column name={Max}}
]
{iter.dat}
\vspace{0.05in}
\caption{Number of iterations used by each method in the tests appearing in Fig.~\ref{fig:err}.}
\label{tab:iter}
\end{table}

\paragraph{Additional tests} 

We performed tests on an additional 44 matrices from the Matrix Function Toolbox~\cite{Higham:MCT}, namely those matrices in the toolbox of size $10 \times 10$ having 2-norm condition number $\kappa_2(A) \le u^{-1}$, where $u = 2^{-53}$ is the unit round-off.  For each matrix $A$, we rescaled $A$ by $e^{i\theta}$ if $A$ had any negative real eigenvalues, with $\theta$ a random number between $0$ and $2\pi$.  

Fig.~\ref{fig:err} shows the relative error $\|\hat{X}-A^{1/2}\|_{\infty} / \|A^{1/2}\|_{\infty}$ committed by each method on the 44 tests, ordered by decreasing condition number $\kappa_{\mathrm{sqrt}}(A)$.  To reduce clutter, the results for the non-Zolotarev iterations (DB, DBp, IN, P-(1,0), P-(4,4), and P-(8,8)) are not plotted individually.  Instead, we identified in each test the smallest and largest relative errors committed among the DB, DBp, IN, P-(1,0), P-(4,4), and P-(8,8) iterations, and plotted these minima and maxima (labelled ``Best iterative'' and ``Worst iterative'' in the legend).  In almost all tests, the Zolotarev iterations achieved relative errors less than or comparable to $u \kappa_{\mathrm{sqrt}}(A)$.  In addition, the Zolotarev iterations tended to produce relative errors closer to the best of the non-Zolotarev iterations than the worst of the non-Zolotarev iterations.

Table~\ref{tab:iter} summarizes the number of iterations used by each method in these tests.  The table reveals that on average, the Zolotarev iteration of type $(m,\ell)$ converged more quickly than the Pad\'e iteration of type $(m,\ell)$ for each $(m,\ell) \in \{(1,0),(4,4),(8,8)\}$.

\section{Conclusion}

We have presented a new family of iterations for computing the matrix square root using recursive constructions of Zolotarev's rational minimax approximants of the square root function.  These iterations are closely related to the Pad\'e iterations, but tend to converge more rapidly, particularly for matrices that have eigenvalues with widely varying magnitudes.  The favorable behavior of the Zolotarev iterations presented here, together with the favorable behavior of their counterparts for the polar decomposition~\cite{nakatsukasa2016computing}, suggests that other matrix functions like the matrix sign function and the matrix $p^{th}$ root may stand to benefit from these types of iterations.

\section*{Acknowledgments}

I wish to thank Yuji Nakatsukasa for introducing me to this topic and for sharing his code for computing the coefficients of Zolotarev's functions.

\bibliography{references}

\begin{thebibliography}{10}

\bibitem{akhiezer1956theory}
{\sc N.~I. Akhiezer}, {\em Theory of Approximation}, Frederick Ungar Publishing
  Corporation, 1956.

\bibitem{akhiezer1990elements}
{\sc N.~I. Akhiezer}, {\em Elements of the Theory of Elliptic Functions},
  vol.~79, American Mathematical Soc., 1990.

\bibitem{beckermann2013optimally}
{\sc B.~Beckermann}, {\em Optimally scaled {N}ewton iterations for the matrix
  square root}, Advances in Matrix Functions and Matrix Equations workshop,
  Manchester, UK, 2013.

\bibitem{beckermann2017singular}
{\sc B.~Beckermann and A.~Townsend}, {\em On the singular values of matrices
  with displacement structure}, SIAM Journal on Matrix Analysis and
  Applications, 38 (2017), pp.~1227--1248.

\bibitem{braess1986nonlinear}
{\sc D.~Braess}, {\em Nonlinear Approximation Theory}, vol.~7, Springer Series
  in Computational Mathematics, 1986.

\bibitem{byers2008new}
{\sc R.~Byers and H.~Xu}, {\em A new scaling for {N}ewton's iteration for the
  polar decomposition and its backward stability}, SIAM Journal on Matrix
  Analysis and Applications, 30 (2008), pp.~822--843.

\bibitem{denman1976matrix}
{\sc E.~D. Denman and A.~N. Beavers~Jr}, {\em The matrix sign function and
  computations in systems}, Applied Mathematics and Computation, 2 (1976),
  pp.~63--94.

\bibitem{gawlik2018backward}
{\sc E.~S. Gawlik, Y.~Nakatsukasa, and B.~D. Sutton}, {\em A backward stable
  algorithm for computing the {CS} decomposition via the polar decomposition},
  (Preprint),  (2018).

\bibitem{guttel2015zolotarev}
{\sc S.~Guttel, E.~Polizzi, P.~T.~P. Tang, and G.~Viaud}, {\em Zolotarev
  quadrature rules and load balancing for the {FEAST} eigensolver}, SIAM
  Journal on Scientific Computing, 37 (2015), pp.~A2100--A2122.

\bibitem{hale2008computing}
{\sc N.~Hale, N.~J. Higham, and L.~N. Trefethen}, {\em Computing {$A^\alpha$},
  {$\log(A)$}, and related matrix functions by contour integrals}, SIAM Journal
  on Numerical Analysis, 46 (2008), pp.~2505--2523.

\bibitem{Higham:MCT}
{\sc N.~J. Higham}, {\em The {Matrix Computation Toolbox}}.
\newblock \path|http://www.ma.man.ac.uk/~higham/mctoolbox|.

\bibitem{higham1986newton}
{\sc N.~J. Higham}, {\em Newton's method for the matrix square root},
  Mathematics of Computation, 46 (1986), pp.~537--549.

\bibitem{higham1997stable}
{\sc N.~J. Higham}, {\em Stable iterations for the matrix square root},
  Numerical Algorithms, 15 (1997), pp.~227--242.

\bibitem{higham2008functions}
{\sc N.~J. Higham}, {\em Functions of Matrices: Theory and Computation}, SIAM,
  2008.

\bibitem{higham2005functions}
{\sc N.~J. Higham, D.~S. Mackey, N.~Mackey, and F.~Tisseur}, {\em Functions
  preserving matrix groups and iterations for the matrix square root}, SIAM
  Journal on Matrix Analysis and Applications, 26 (2005), pp.~849--877.

\bibitem{hogben2016handbook}
{\sc L.~Hogben}, {\em Handbook of Linear Algebra}, CRC Press, 2016.

\bibitem{iannazzo2003note}
{\sc B.~Iannazzo}, {\em A note on computing the matrix square root}, Calcolo,
  40 (2003), pp.~273--283.

\bibitem{kressner2017fast}
{\sc D.~Kressner and A.~Susnjara}, {\em Fast computation of spectral projectors
  of banded matrices}, SIAM Journal on Matrix Analysis and Applications, 38
  (2017), pp.~984--1009.

\bibitem{bailly2000optimal}
{\sc B.~Le~Bailly and J.~Thiran}, {\em Optimal rational functions for the
  generalized {Z}olotarev problem in the complex plane}, SIAM Journal on
  Numerical Analysis, 38 (2000), pp.~1409--1424.

\bibitem{li2017spectrum}
{\sc Y.~Li and H.~Yang}, {\em Spectrum slicing for sparse {H}ermitian definite
  matrices based on {Z}olotarev's functions}, arXiv preprint arXiv:1701.08935,
  (2017).

\bibitem{meini2004matrix}
{\sc B.~Meini}, {\em The matrix square root from a new functional perspective:
  theoretical results and computational issues}, SIAM Journal on Matrix
  Analysis and Applications, 26 (2004), pp.~362--376.

\bibitem{nakatsukasa2016computing}
{\sc Y.~Nakatsukasa and R.~W. Freund}, {\em Computing fundamental matrix
  decompositions accurately via the matrix sign function in two iterations:
  {T}he power of {Z}olotarev's functions}, SIAM Review, 58 (2016),
  pp.~461--493.

\bibitem{ninomiya1970best}
{\sc I.~Ninomiya}, {\em Best rational starting approximations and improved
  {N}ewton iteration for the square root}, Mathematics of Computation, 24
  (1970), pp.~391--404.

\bibitem{rutishauser1963betrachtungen}
{\sc H.~Rutishauser}, {\em Betrachtungen zur quadratwurzeliteration},
  Monatshefte f{\"u}r Mathematik, 67 (1963), pp.~452--464.

\bibitem{trefethen1985convergence}
{\sc L.~N. Trefethen and M.~H. Gutknecht}, {\em On convergence and degeneracy
  in rational {P}ad{\'e} and {C}hebyshev approximation}, SIAM Journal on
  Mathematical Analysis, 16 (1985), pp.~198--210.

\bibitem{wachspress2013adi}
{\sc E.~Wachspress}, {\em The {ADI} Model Problem}, Springer, 2013.

\bibitem{zolotarev1877applications}
{\sc E.~I. Zolotarev}, {\em Applications of elliptic functions to problems of
  functions deviating least and most from zero}, Zapiski St-Petersburg Akad.
  Nauk, 30 (1877), pp.~1--59.

\end{thebibliography}

\end{document}